\documentclass[reqno] {amsart}
\usepackage{amsmath}
\usepackage{amscd}
\usepackage{amsthm}
\usepackage{amssymb}
\usepackage{latexsym}
\usepackage{enumerate}

\newtheorem{thm}{Theorem}[section]
\newtheorem{prop}{Proposition}[section]
\newtheorem{lem}{Lemma}[section]

\newtheorem{rem}{Remark}[section]

\title{The radial curvature of an end that makes eigenvalues vanish in the essential spectrum I}

\author{Hironori Kumura}

\date{}

\begin{document}

\maketitle

\begin{abstract}
The concern of this paper is to clarify a relationship between the curvatures at infinity and the spectral structure of the Laplacian. 
In particular, this paper discusses the question of whether there is an eigenvalue of the Laplacian embedded in the essential spectrum or not. 
The borderline-behavior of the radial curvatures for this problem will be determined: we will assume that the radial curvature $K_{\rm rad.}$ of an end converges to a constant $-1$ at infinity with the decay order $K_{\rm rad.} + 1 = o(r^{-1})$ and prove the absence of eigenvalues embedded in the essential spectrum. Furthermore, in order to show that this decay order $K_{\rm rad.} + 1 = o( r^{-1} )$ is sharp, we will construct a manifold with the radial curvature decay $K_{\rm rad.} + 1 = O(r^{-1})$ and with an eigenvalue $\frac{(n-1)^2}{4} + 1$ embedded in the essential spectrum $[ \frac{(n-1)^2}{4} , \infty )$ of the Laplacian. 
\end{abstract}
\section{Introduction}

The Laplace-Beltrami operator $\Delta$ on a noncompact complete Riemannian manifold $(M,g)$ is essentially self-adjoint on $C^{\infty}_0(M)$; the spectral structure of its self-adjoint extension to $L^2(M)$ and the curvatures of the manifold $M$ are closely related to each other and their relationship has been studied by several authors. 
Especially, the problem of the absence of eigenvalues was discussed in [2,3,4,5,6,7,9,10,13,16, \linebreak 
18,20]. 

This paper will treat the case that the curvature $K$ of $(M,g)$ converges to a constant $-1$ at infinity. 
The other case that $K$ converges to $0$ at infinity will be treated in the sequel. 
We will begin by recalling the previous works and their decay conditions imposed on $K+1$ which ensure the absence of eigenvalues greater than $\frac{(n-1)^2}{4}$. 
In the case that $(M,g)$ is rotationally symmetric, the condition imposed in Pinsky \cite{P} is that $\mathrm{dim} \,M = 2$, $ K_{\rm rad.} \le 0$ on $(M,g)$, $K_{\rm rad.}\le -1~(r\ge r_0)$, and $\int_{r_0}^{\infty} | K_{\rm rad.}+1|\,dr<\infty$, where $K_{\rm rad.}$ is the radial curvature (see \cite{G-W} ) with respect to the origin. 
When the metric is not necessarily rotationally symmetric, the condition imposed in Donnelly \cite{D2} is that $(M,g)$ is a simply connected negatively curved manifold, $\int_1^{\infty} r^{\beta} |K_{\rm rad.} + 1| \,dr < \infty$ and $ \lim_{r\to \infty} r^{\beta} |K_{\rm rad.}+1| = 0$, where $r$ denotes the distance to an arbitrarily fixed point $p$, $K_{\rm rad.}$ stands for the radial curvature with respect to $p$, and $\beta>2$ is a constant. 
Roughly speaking, Donnelly's curvature condition is $K_{\rm rad.} + 1 = O( r^{-3- \varepsilon} )$. 
In this paper, we will examine the growth property of solutions to eigenvalue equation and determine the borderline-behavior of the radial curvatures that makes the eigenvalues in the essential spectrum vanish. 
Indeed, as a special case, we will prove the absence of eigenvalues greater than $\frac{(n-1)^2}{4}$ under the assumption that $K_{\rm rad.} + 1 = o(r^{-1})$; this paper also construct an example which shows that this curvature condition $K_{\rm rad.} + 1 = o(r^{-1})$ is sharp. 
Thus, the borderline-behavior of the radial curvatures that makes the eigenvalues embedded in the essential spectrum will be seen to be $o(r^{-1})$. 

We state our results more precisely: let $(M,g)$ be an $n$-dimensional noncompact complete Riemannian manifold and $U$ an open subset of $M$. 
We shall say that $M-U$ is an {\it end with radial coordinates} if and only if the boundary $\partial U$ is compact, connected, and $C^{\infty}$ and the outward normal exponential map $\exp_{\partial U}^{\perp} : N^{+}(\partial U) \to M - \overline{U}$ induces a diffeomorphism, where $N^{+}(\partial U) = \{ v\in T(\partial U) \mid v {\rm ~is~outward~normal~to~}\partial U \}$. 
Note that $U$ is not necessarily relatively compact. 
Let $r$ denote the distance function from $\partial U$ defined on the end $M-\overline{U}$. 
We shall say that a $2$-plane $\pi \subset T_xM$ $(x\in M-\overline{U})$ is {\it radial} if $\pi$ contains $\nabla r$, and, by the {\it radial curvature}, we mean the restriction of the sectional curvature to all the radial planes. 
In the sequel, the following notations will be used:
\begin{align*}
   & B(s,t) = \{ x\in M - \overline{U} \mid s < r(x) < t \} \quad 
              \mathrm{for}~~0 \le s < t ;\\
   & B(s,\infty) = \{x \in M-\overline{U}~|~s<r(x)\} \quad 
              \mathrm{for}~~0 \le s <\infty;\\
   & S(t) = \{x\in M-\overline{U}~|~r(x)=t\} \quad 
              \mathrm{for}~~0 \le t<\infty ;\\
   & \sigma(-\Delta) =  \,{\rm the~spectrum~of}~-\Delta;\\
   & \sigma_{\rm p}(-\Delta) =  \,{\rm the~set~of~all~eigenvalues~of}~
             - \Delta;\\
   & \sigma_{\rm ess}(-\Delta) = \,{\rm the~essential~spectrum~of}~- \Delta;\\
   & K_{{\rm rad.}} =  \,{\rm the~radial~curvature~on~}M-U.
\end{align*}
Moreover, we denote the Riemannian measure of $(M,g)$ by $dv_g$, and the measure on each $S(t)~(t>0)$ induced from $dv_g$ simply by $dA$. 

In this paper, we shall consider the eigenvalue equation 
$$
     \Delta f + \alpha f = 0
$$
on an end $M-\overline{U}$ and drive a growth estimate at infinity of solutions $f$, from which will follow the absence of eigenvalues in the essential spectrum. 

First, we shall state our theorem in terms of the shape operators $\nabla dr$ of the level hypersurfaces $\{S(r)\}_{r\ge 0}$ and the lower bound of the Ricci curvature of the radial direction; recall that the shape operators $\nabla dr$ of $\{S(r)\}_{r\ge 0}$ describes the metric-growth on $\{S(r)\}_{r\ge 0}$:
\begin{thm}
Let $(M,g)$ be an $n$-dimensional noncompact complete Riemannian manifold and $U$ an open subset of $M$. 
Assume that the end $E:= M - U $ has radial coordinates and write $r={\rm dist}\,(U,*)$. 
We also assume that there exists a constant $ r_0 > 0 $ such that 
\begin{align}
   \left( 1-\frac{A_1}{r} \right) \widetilde{g}
  \le 
   \nabla dr \le \left( 1+\frac{B_1}{r}\right) \widetilde{g}
   \quad {\rm on}~~B(r_0,\infty)  \tag{$*_1$}
\end{align}
and 
\begin{align}
  {\rm Ric}\,(\nabla r, \nabla r) 
  \ge - (n-1) \left( 1 + \frac{b_1}{r} \right) 
  \quad {\rm on}~~B(r_0,\infty) , \tag{$*_2$}
\end{align}
where we set $\widetilde{g}= g - dr\otimes dr$ for simplicity, and $A_1$, $B_1$, and $b_1$ are positive constants. 
Let $\alpha > \frac{(n-1)^2}{4}$ be a constant and $f$ a not identically vanishing solution to the equation$\,:$
\begin{align*}
   \Delta f + \alpha f = 0  \quad\qquad {\rm on}~~B(r_0,\infty).
\end{align*}
Assume that constants $\gamma >0$, $A_1$, $B_1$, and $b_1$ satisfy 
\begin{align}
  & 1 - \widehat{A}_1 > 0 ; 
    \quad 2 \gamma > \widehat{A}_1 + \widehat{B}_1 ; \tag{$1.1$} \\
  & \alpha - \frac{(n-1)^2}{4} 
    >
   (n-1) \left( 2 \widehat{A}_1 + \widehat{b}_1 \right) 
    \cdot m_1( \gamma, A_1, B_1 ),  \tag{$1.2$}
\end{align}
where we set $\widehat{A}_1 = (n-1) A_1$, $\widehat{B}_1 = (n-1) B_1$, 
$\widehat{b}_1 = (n-1) b_1$, and 
\begin{align*}
   m_1( \gamma, A_1, B_1 ) = 
   \max \left\{ \frac{1}{ 2( 1 - \widehat{A}_1 )}, 
         \frac{1}{ 2 \gamma - \widehat{A}_1 - \widehat{B}_1 } \right\} .
\end{align*}
Then, we have
\begin{equation*}
   \liminf_{t\to \infty}~t^{\gamma}\int_{S(t)}
   \left\{ \left( \frac{\partial f}{\partial r} \right) ^2 
   + |f|^2 \right\} \,dA  = \infty.
\end{equation*}
\end{thm}
Next, we shall state our theorem in terms of radial curvatures: the conditions $(1.3)$ and $(1.4)$ in the following theorem implies $(*_1)$ and $(*_2)$ with $ b_1 = 2B_1 $ in Theorem $1.1$ except for unnecessary lower order terms (Proposition $2.1$). 
Hence, the following theorem follows from Theorem $1.1$:
\begin{thm}
Let $(M,g)$ be an $n$-dimensional noncompact complete Riemannian manifold and $U$ an open subset of $M$. 
Assume that the end $E := M - U$ has radial coordinates and write $r={\rm dist}\,(U,*)$. 
We also assume that there exists a constant $ r_0 > 0 $ such that 
\begin{align}
   \nabla dr \ge 0 \qquad {\rm on}~~S(r_0)  \tag{$1.3$}
\end{align}
and 
\begin{align}
  - 1 - \frac{ 2B_1 }{r} \le K_{{\rm rad.}} \le - 1 + \frac{ 2 A_1 }{r}
  \le 0
  \qquad {\rm on}~~B(r_0,\infty) . \tag{$1.4$}
\end{align}
Let $\alpha > \frac{(n-1)^2}{4}$ be a constant and $f$ a not identically vanishing solution to the equation$\,:$
\begin{align*}
   \Delta f + \alpha f = 0  \quad\qquad {\rm on}~~B(r_0,\infty).
\end{align*}
Assume that constants $\gamma >0$, $A_1$, $B_1$, $b_1$ satisfy 
\begin{align*}
  & 1 - \widehat{A}_1 > 0 ; 
    \quad 2 \gamma > \widehat{A}_1 + \widehat{B}_1 ; \tag{$1.1$} \\
  & \alpha - \frac{(n-1)^2}{4} 
    >
    2 (n-1) \left( \widehat{A}_1 + \widehat{B}_1 \right) 
    \cdot m_1( \gamma, A_1, B_1 ).  \tag{$1.2$}
\end{align*}
Then, we have
\begin{equation*}
   \liminf_{t\to \infty}~t^{\gamma}\int_{S(t)}
   \left\{ \left( \frac{\partial f}{\partial r} \right) ^2 
   + |f|^2 \right\} \,dA  = \infty.
\end{equation*}
\end{thm}
In the theorems above, if $A_1$ and $B_1$ converge to $0$, the conditions $(1.1)$ and $(1.2)$ are satisfied for any $\gamma >0$; thus, we get the following two theorems by taking $r_0$ successively large:
\begin{thm}
Let $(M,g)$ be an $n$-dimensional noncompact complete Riemannian manifold and $U$ an open subset of $M$. 
Assume that the end $E:= M - U $ has radial coordinates and write $r={\rm dist}\,(U,*)$. 
We also assume that there exists a constant $ r_0 > 0 $ such that 
\begin{align*}
   \left( 1-\frac{A(r)}{r} \right) \widetilde{g}
  \le 
   \nabla dr 
  \le \left( 1+\frac{B(r)}{r}\right) \widetilde{g}
   \quad {\rm on}~~B(r_0,\infty)  
\end{align*}
and 
\begin{align*}
  {\rm Ric}\,(\nabla r, \nabla r) 
  \ge - (n-1) \left( 1 + \frac{b(r)}{r} \right) 
  \quad {\rm on}~~B(r_0,\infty) , 
\end{align*}
where we set $\widetilde{g}= g - dr\otimes dr$ and $A(r)$, $B(r)$ and $b(r)$ are positive-valued continuous function of $r \in [r_0, \infty)$ satisfying $\lim_{r \to \infty} A(r) = \lim_{r \to \infty} B(r) = \lim_{r \to \infty} b(r) = 0$. 
Let $\alpha > \frac{(n-1)^2}{4}$ be a constant and $f$ a not identically vanishing solution to the equation$\,:$
\begin{align*}
   \Delta f + \alpha f = 0  \quad\qquad {\rm on}~~B(r_0,\infty).
\end{align*}
Then, we have for any $\gamma >0$
\begin{equation*}
   \liminf_{t\to \infty}~t^{\gamma}\int_{S(t)}
   \left\{ \left( \frac{\partial f}{\partial r} \right) ^2 
   + |f|^2 \right\} \,dA  = \infty.
\end{equation*}
\end{thm}
\begin{thm}
Let $(M,g)$ be an $n$-dimensional noncompact complete Riemannian manifold and $U$ an open subset of $M$. 
Assume that the end $E : = M - U$ has radial coordinates. 
Moreover, we assume that the end $E$ satisfies the following conditions$\,:$ there exists $r_0 > 0$ such that 
\begin{align}
   & \nabla dr \ge 0 \qquad \mathrm{on}~~S(r_0);  \tag{$1.3$} \\
   -1-\frac{b(r)}{r} \le \,\,& K_{{\rm rad.}} \le -1+\frac{a(r)}{r} \le 0  \qquad \mathrm{on}~~ B(r_0,\infty ),  \tag{$1.5$}
\end{align}
where $a(r)$ and $b(r)$ are positive-valued continuous functions of $r \in [r_0,\infty)$ satisfying $ \lim_{r\to \infty} a(r) = \lim_{r\to \infty} b(r) = 0$. 
If $\alpha >\frac{(n-1)^2}{4}$ and if $f$ is a not identically vanishing solution to $\Delta f + \alpha f = 0$ on $E$, then we have for any $\gamma >0$
\begin{equation*}
   \liminf_{t\to \infty}~t^{\gamma}\int_{S(t)}
   \left\{ \left( \frac{\partial f}{\partial r} \right)^2 + |f|^2 \right\} 
   \,dA = \infty.
\end{equation*}
\end{thm}
Let us now look at our theorems from the viewpoint of the spectral structure; by the growth property of solutions of eigenvalue equation mentioned in theorems above, we get the following theorems:
\begin{thm}
Let $(M,g)$ be an $n$-dimensional complete Riemannian manifold. 
Assume that $(M,g)$ has at least one end mentioned either in Theorem $1.1$ with $\gamma = 1$ or in Theorem $1.2$ with $\gamma = 1$. 
Then $\sigma_{{\rm ess}}(-\Delta) \supseteq \left[\frac{(n-1)^2}{4},\infty\right)$ and any 
\begin{align*}
     \alpha 
   > \frac{(n-1)^2}{4} + 
     (n-1) \left( 2 \widehat{A}_1 + \widehat{b}_1 \right)
   \cdot m_1( 1 , A_1, B_1 )
\end{align*}
is not eigenvalue of $-\Delta $. 
\end{thm}
\begin{thm}
Let $(M,g)$ be an $n$-dimensional complete Riemannian manifold. 
Assume that $(M,g)$ has at least one end mentioned either in Theorem $1.3$ or in Theorem $1.4$. 
Then $ \sigma_{{\rm ess}} (-\Delta) \supseteq \left[ \frac{(n-1)^2}{4}, \infty \right) $ and any $ \alpha > \frac{(n-1)^2}{4}$ is not eigenvalue of $-\Delta $. \end{thm}
\begin{thm}
Let $(M,g)$ be an $n$-dimensional complete Riemannian manifold and $U$ a relatively compact open subset of $M$ with $C^{\infty}$-boundary. 
We assume that the components of $M-\overline{U}$ consists of the disjoint union of finite number of ends $E_1,E_2,\cdots,E_l$ with radial coordinates. 
Assume that each end $E_i~(1\le i \le l)$ satisfies 
\begin{align*}
   &\nabla dr\ge 0 \qquad \mathrm{on}~~\partial E_i;\\
   -k_i-\frac{b_i(r)}{r}
  \le \,\,
   & K_{{\rm rad.}}\le -k_i + \frac{a_i(r)}{r}\le 0 \qquad \mathrm{on}~~ E_i,
\end{align*}
where $r$ is the distance function to $U$, $k_i>0$ is a constant, and, $a_i(r)$ and $b_i(r)$ are positive continuous function of $r$ satisfying $\lim_{r\to \infty} a_i(r) = \lim_{r\to \infty} b_i(r) = 0$. 
Then, $\sigma_{{\rm ess}}(-\Delta) = \left[ \frac{(n-1)^2 k_{\min}}{4}, \infty\right)$, and any $\alpha > \frac{(n-1)^2 k_{\min}}{4}$ is not eigenvalue of $- \Delta $, where we set $k_{\min} = \min \{ k_i \mid 1 \le i \le l \}$. 
\end{thm}
In section $8$, we shall construct a rotationally symmetric manifold which has eigenvalue $\frac{(n-1)^2}{4} + 1$ in the essential spectrum $\left[ \frac{(n-1)^2}{4},\infty \right)$ and has the radial curvature decay $ K_{\rm rad.} + 1 = O(r^{-1})$: 
\begin{thm}
There exists a rotationally symmetric manifold 
$(M,g) = \bigl( {\bf R}^n, dr^2 + f^2(r) g_{S^{n-1}(1)} \bigr)$ with the following three properties$\,:$
\begin{enumerate}[$(1)$]
   \item $ \lim_{r\to \infty} \left| \nabla dr-(g-dr\otimes dr) \right| = 0$, and hence, $\sigma_{{\rm ess}}(-\Delta)=\left[\frac{(n-1)^2}{4},\infty \right)$$;$   \item $\sigma_{\rm p} ( -\Delta ) \cap \left( \frac{(n-1)^2}{4} , \infty \right)= \left\{ \left(\frac{(n-1)^2}{4} \right) + 1 \right\}$$;$
   \item $ K_{{\rm rad.}}+1=O(r^{-1})$ as $r\to \infty$.
\end{enumerate}
\end{thm}
Theorem $1.8$ shows that the curvature decay condition $ K_{\rm rad.} + 1 = o(r^{-1})$ in Theorem $1.4$ and Theorem $1.6$ is sharp. 
Thus, the borderline-behavior of the radial curvatures that makes the eigenvalues in the essential spectrum vanish is $o(r^{-1})$. 

Hessian comparison theorem (see Kasue \cite{Kasue}) in Riemannian geometry is an important ingredients of our proof, and our method is a modification of solutions of Kato \cite{Kato}, Eidus \cite{E}, Roze \cite{R} and Mochizuki \cite{M} to the analogous problem for the Schr\"odinger equation on Euclidian space. 

\vspace{3mm}
The author would like to express his gratitude to Professor Minoru Murata; he kindly informed the author of several facts about the analogous results for the Schr\"odinger equation on Euclidian space. 

\section{Geometric situation}

In this section, we shall confirm our geometric situation. 

The Hessian comparison theorem (see Kasue \cite{Kasue}) implies the following:
\begin{prop}
Let $(M,g)$ be an $n$-dimensional noncompact complete Riemannian manifold and $U$ an open subset of $M$. 
Assume that the end $E := M-U$ has radial coordinates and set $r (x) = {\rm dist}\,(U,x)$ for $x \in E$. 
We assume that there exist constants $r_0>0$, $A_1>0$, and $B_1>0$ such that 
\begin{align}
  & \nabla dr \ge 0  \qquad \mathrm{on} ~~S(r_0); \tag{$1.3$} \\
   - 1 - \frac{2B_1}{r}\le \,\,
  & K_{{\rm rad.}}
    \le - 1 + \frac{2A_1}{r} \le 0  
    \qquad \mathrm{on} ~~B(r_0,\infty ) \tag{$1.4$}.
\end{align}
Then we have 
\begin{align*}
   \left\{ 1 - \frac{A_1}{r} + O \left( \frac{1}{r^3}\right ) \right\} 
   \widetilde{g}
   \le 
   \nabla dr 
   \le 
   \left\{ 1 + \frac{B_1}{r} + O \left( \frac{1}{r^3}\right ) \right\} 
   \widetilde{g}
\end{align*}
on $B(r_0,\infty )$, where we set $\widetilde{g} = g - dr \otimes dr$, for simplicity. 
\end{prop}
\begin{proof}
Let $\widehat{{\mathcal A}}$ denote the shape operators of the level hypersurfaces $\{ S(r) \}_{ r\ge r_0 }$ with respect to the inward unit normal $ - \nabla r: = - {\rm grad}~r$; 
that is, $\langle \widehat{{\mathcal A}}_x u,v \rangle = \nabla dr (u,v)$ for $x\in E$, $u,v \in T_x S ( r(x) ) $. 
Recall that $\widehat{{\mathcal A}}$ satisfies the Riccati-type equation along each normal geodesic $[0,\infty) \ni t \to \exp_{\partial U}^{\perp} (tu)$ $(u \in N^{+}(\partial U),~|u|=1)$:
\begin{align}
   \nabla_{\nabla r} \widehat{{\mathcal A}} + \widehat{{\mathcal A}}\,^2 + 
   R(*,\nabla r)\nabla r = 0  \qquad {\rm on}~~\nabla r^{\perp} , 
\end{align}
where $\nabla r^{\perp} := TS = \{u \in T_xS(r(x)) \mid x \in M-U \}$ and $R$ stands for the Riemannian curvature tensor: 
$R(X.Y)Z=\nabla_X\nabla_Y Z - \nabla_Y\nabla_X Z - \nabla_{[X,Y]}Z$. 
Note that 
$$
   \langle R(u,\nabla r)\nabla r ,u \rangle 
   = K_{{\rm rad.}}(\{\nabla r,u\}_{{\bf R}}),
$$ 
where $\{\nabla r,u\}_{{\bf R}}$ is the $2$-plane spanned by $\nabla r$ and $u\in \nabla r^{\perp}$. 
Set 
\begin{align}
 K_1(r) := - 1 + \frac{2A_1}{r}\,(\le 0);\quad K_2(r) = - 1 - \frac{2B_1}{r} 
\end{align}
and consider solutions $f_1(r)$ and $f_2(r)$ to the ordinary differential equations
\begin{align}
  f'(r) + f(r)^2 + K_i(r) = 0 \qquad (i=1,2) 
\end{align}
with the initial conditions
\begin{align}
    f_1(r_0)=0;~~f_2(r_0)
    =\max\left\{ 
    \frac{\langle \widehat{{\mathcal A}}_x u,u \rangle}{\langle u, u \rangle} 
    \, \bigg| \, x\in S(r_0), 0 \neq u\in T_xS(r_0) \right\}. 
\end{align}
With these understanding, the radial-curvature-assumption $(1.4)$ and initial conditions $(1.3)$ and $(4)$ are written by $f_1(r_0)=0 \le \widehat{{\mathcal A}} \,\big|_{S(r_0)} \le f_2(r_0)$ and $K_2(r) \le \, K_{{\rm rad.}}\le K_1(r) \le 0$ on  $B(r_0,\infty)$, and hence, applying the Hessian comparison theorem (see Kasue \cite{Kasue}) to the equations $(1)$ and $(3)$, we see that the eigenvalues of the symmetric operators $\widehat{{\mathcal A}}$ on $\nabla r^{\perp}$ is pinched between two numbers $f_1(r)$ and $f_2(r)$:
\begin{align}
  f_1(r) \le \widehat{{\mathcal A}} \le f_2(r) 
  \qquad {\rm on}~~B(r_0,\infty) .
\end{align}
Since $f_1(r_0) = 0$ and $K_1(r) \le 0$, the comparison theorem implies that 
\begin{align}
    f_1(r) \ge 0 \qquad {\rm on}~~ [r_0,\infty)  . 
\end{align}
In view of $(2)$ and $(6)$, we see that solutions $f_1$ and $f_2$ to $(3)$ have the following asymptotic behavior:
\begin{align}
 f_1(r) = 1 - \frac{A_1}{r} + O \left(\frac{1}{r^3}\right);~~
 f_2(r) = 1 + \frac{B_1}{r} + O \left(\frac{1}{r^3}\right). 
\end{align}
Proposition $2.1$ follows from $(5)$ and $(7)$. 
\end{proof}
Proposition $2.1$ shows that Theorem $1.1$ implies Theorem $1.2$. 

We can prove the following proposition in the same way:
\begin{prop}
In Proposition $2.1$, replace the condition $(1.4)$ with the following:
\begin{align}
  & - 1 - \frac{b(r)}{r}\le \,\,
    K_{{\rm rad.}}
    \le - 1 + \frac{a(r)}{r} \le 0  
    \qquad \mathrm{on} ~~B(r_0,\infty ) ; \tag{$1.5$} \\
  & \hspace{10mm} 
    \lim_{r\to \infty} a(r) =  \lim_{r\to \infty} b(r) = 0. \nonumber
\end{align}
Then we have
\begin{align*}
  & \left( 1 - \frac{A(r)}{r} \right) \widetilde{g} 
    \le \nabla dr \le 
    \left( 1 + \frac{B(r)}{r} \right) \widetilde{g}
    \qquad \mathrm{on} ~~B(r_0,\infty ); \\
  & \hspace{10mm}
    \lim_{r\to \infty} A(r) =  \lim_{r\to \infty} B(r) = 0.
\end{align*}
\end{prop}
If $\Delta r$ converges to a constant $n-1$ as $r$ tends to infinity on one end, then this end produces the essential spectrum $ \left[ \frac{(n-1)^2}{4} , \infty \right) \subseteq \sigma_{{\rm ess}}(-\Delta ) $ (see \cite{K1}). 

The following identity will play an important role in our proof:
\begin{prop}
Let $(M,g)$ be an $n$-dimensional Riemannian manifold and $P$ a closed subset of $M$. 
We assume that the distance function $r = {\rm dist}\, (P,*)$ is $C^{\infty}$ on an open neighborhood $W$ of $P$. 
Then, on $W$, we have
\begin{align}
   -\frac{\partial (\Delta r)}{\partial r} 
   =  |\nabla dr|^2 + \mathrm{Ric}\,(\nabla r,\nabla r)  \tag{$2.1$}
\end{align}
\end{prop}
\begin{rem}
Although the identity $(2.1)$ is proved by substituting $u_1 = r$ in Weitzenb\"ock formula$:$
\begin{align*}
  \frac{1}{2} \Delta \bigl( |\nabla u_1 |^2 \bigr) 
  = |\nabla d u_1|^2 + \langle \nabla u_1 , \nabla \Delta u_1 \rangle 
  + {\rm Ric} \,( \nabla u_1, \nabla u_1),
\end{align*}
 the readers should regard this identity $(2.10)$ as the trace of the Riccati-type equation $(1)$. 
We should notice that this important identity is used in the proof of Bishop-Gromov comparison theorem. 
\end{rem}

\section{Analytic propositions}

In this section, we shall prepare some analytic propositions for the proof of Theorem $1.1$. 

The purpose of this study is to examine solutions to the eigenvalue equation 
\begin{align*}
    \Delta f + \alpha f = 0 \quad\qquad {\rm on}~~E := M-\overline{U},
\end{align*}
where $\alpha $ is a constant satisfying $\alpha  >\frac{(n-1)^2}{4}$. 
In the sequel, we shall derive a growth estimate at infinity of the solution $f$, from which the absence of the eigenvalue will follow. 

Let us set $c = \frac{(n-1)}{2}$ and transform the operator $\Delta+c^2$ and Riemannian measure $dv_g$ into the new operator $L = e^{cr} (\Delta + c^2) e^{-cr} = \Delta - 2c \frac{\partial}{\partial r} + c( 2c - \Delta r )$ and new measure $e^{-2cr} dv_g$, respectively: 
\begin{align*}
  \begin{CD}
    L^2(E,dv_g)           @>{ -(\Delta + c ^2)}>>       L^2(E,dv_g)          \\
    @V{e^{cr}}VV                         @VV{e^{cr}}V  \\
    L^2(E,e^{-2cr}dv_g)   @>>{ - L }>           L^2(E,e^{-2cr}dv_g)
  \end{CD}
\end{align*}
Note that the multiplying operator $e^{cr}:L^2(E,dv_g)\ni h \mapsto e^{cr}h \in L^2(E,e^{-2cr}dv_g)$ is a unitary operator. 
Thus, in order to show that $-\Delta$ has no nontrivial $L^2(E,dv_g)$-eigenfunction with eigenvalue greater than $\frac{(n-1)^2}{4}$, it will suffice to show that $-L$ has no nontrivial $L^2(E,e^{-2cr}dv_g)$-eigenfunction with positive eigenvalue. 
For simplicity, we put the new measures as follows:
\begin{align*}
  d\mu_c := e^{-2cr} dv_g ; ~~dA_c := e^{-2cr}dA,
\end{align*}
where recall that $dA$ is the measure on each level surface $S(t)~~(t>0)$ induced from the Riemannian measure $dv_g$ 

Now, let $\lambda$ be a positive constant and $u$ a solution to the following equation:
\begin{equation}
  Lu + \lambda u=0  \quad\qquad {\rm on}~~E,
\end{equation}
that is, 
\begin{align}
   \Delta u - 2c \frac{\partial u}{\partial r} + c(2c-\Delta r)u + \lambda u
   = 0 \quad\qquad {\rm on}~~E. \tag{8'}
\end{align}
Let $\rho (r)$ be a $C^{\infty}$ function of $r \in [r_0,\infty)$, and put 
\begin{align*}
  v(x) = \exp \bigl(\rho (r(x)) \bigr)u(x) \qquad {\rm for}~~x \in E.
\end{align*}
Since $\Delta u=e^{-\rho}\left\{ \Delta v - 2\rho'(r) \frac{\partial v}{\partial r} + \left(|\rho'(r)|^2 - \Delta \rho\right)v\right\}$ and $\frac{\partial u}{\partial r}=e^{-\rho}\left\{ \frac{\partial v}{\partial r} - \rho'(r)v\right\}$, substituting these equations to $(8$'), we see that $v$ satisfies the following  equation on $B(r_0,\infty)$:
\begin{align}
    & \Delta v - 2 \bigl( \rho '(r) + c \bigr) \frac{\partial v}{\partial r} 
     + qv = 0,\\
    & q = |\nabla \rho |^2 - \Delta \rho + \lambda + 
      c\left( 2\rho '(r) + 2c - \Delta r \right) \nonumber\\
    & \hspace{2mm} = |\rho'(r)|^2 - \rho''(r) + 
    \left( 2c - \Delta r \right) \left( \rho'(r) + c \right) + \lambda,
\end{align}
where we set $\nabla v={\rm grad}\,v$. 
In order to prove Theorem $1.1$, we will prepare three Propositions.  
The first is the integration-by-parts-lemma: 
\begin{lem}[integration by parts]
For any $f_1, g_1 \in C^{\infty}(M-\overline{U})$, $C^{\infty}$-vector field $X$, and $r_0< s< t$, we have
\begin{align*}
   &\int_{B(s,t)}\left(\Delta f_1 \right) h_1 \,d\mu_c\\
  =&\left( \int_{S(t)} - \int_{S(s)}\right)
   \frac{ \partial f_1 }{\partial r} h_1 \, dA_c 
   - \int_{B(s,t)}
   \langle \nabla f_1, \nabla h_1 - 2c\, h_1 \nabla r \rangle \,d\mu_c
\end{align*}
and 
\begin{align*}
   \int_{B(s,t)} ({\rm div}\, X)\, d\mu_c
  =\left( \int_{S(t)} - \int_{S(s)}\right)
    \langle X,\nabla r \rangle \, dA_c 
    + 2c \int_{B(s,t)} \langle X,\nabla r \rangle \,d\mu_c.
\end{align*}
\end{lem}
\begin{prop}
For any $\psi \in C^{\infty}(M-\overline{U})$ and $r_0< s< t$, we have
\begin{align*}
    &  \int_{B(s,t)}\left\{ |\nabla v|^2 - q |v|^2 \right\}\psi \,d\mu_c  \\
  = &  \left( \int_{S(t)}-\int_{S(s)}\right )
      \frac{\partial v}{\partial r} \psi v \,dA_c - \int_{B(s,t)}
      \left\langle 
      \nabla \psi + 2 \psi \rho '(r) \nabla r , \nabla v 
      \right\rangle v 
      \,d\mu_c.
\end{align*}
\end{prop}
\begin{proof}
Multiply the equation $(9)$ by $\psi v$ and integrate it over $B(s,t)$ with respect to the measure $d\mu_c = e^{-2cr} \,dv_g$.
Then, the integration-by-parts-lemma with $f_1 = v$ and $g_1 = \psi v$ yields Proposition $3.1$.
\end{proof}
\begin{prop}
For any $r_0< s< t$ and $\gamma \in \mathbf{R}$, we have
\begin{align}
   & \left( \int_{S(t)} - \int_{S(s)} \right) r^{\gamma}
     \left\{ \left(\frac{\partial v}{\partial r} \right)^2
     - \frac{1}{2}|\nabla v|^2 + \frac{1}{2}q|v|^2 \right\} \,dA_c \nonumber \\
 = & \int_{B(s,t)} r^{\gamma -1} \left\{ r(\nabla dr)(\nabla v,\nabla v)
    - \frac{1}{2} ( \gamma + r\Delta r - 2cr ) \left( |\nabla v|^2 
    - \left( \frac{\partial v}{\partial r} \right)^2 \right) \right\} 
     \,d\mu _c  \nonumber \\
   & + \int_{B(s,t)}r^{\gamma -1} 
     \left\{ \frac{1}{2}( \gamma - r\Delta r + 2cr ) + 2r\rho'(r) \right\}
     \left( \frac{\partial v}{\partial r} \right)^2 \,d\mu _c  \nonumber \\
   &  + \frac{1}{2} \int_{B(s,t)}r^{\gamma -1} 
     \left\{ 
     ( \gamma + r\Delta r - 2cr ) q + r \frac{\partial q}{\partial r} 
     \right\}|v|^2 \,d\mu _c .
\end{align}
\end{prop}
\begin{proof}
Let us multiply the equation $(9)$ by $\frac{\partial v}{\partial r}$. 
Then, three identities 
\begin{align*}
  {\rm div} \left( \frac{\partial v}{\partial r} \nabla v \right)
  & = \frac{\partial v}{\partial r} \Delta v 
      + \left\langle 
        \nabla v,\nabla \left(\frac{\partial v}{\partial r}\right) 
        \right\rangle ;\\
  \left\langle \nabla v,
     \nabla \left( \frac{\partial v}{\partial r} \right) \right\rangle 
  & = (\nabla v) \left\langle \nabla r, \nabla v \right\rangle 
     =(\nabla dr) (\nabla v, \nabla v) + (\nabla dv) (\nabla r, \nabla v);\\
  (\nabla dv)(\nabla r,\nabla v)
  & = \frac{1}{2} {\rm div}\,\left( |\nabla v|^2 \nabla r \right)
      - \frac{1}{2} |\nabla v|^2 \Delta r,
\end{align*}
yield 
\begin{align}
    - \mathrm{div} \left( \frac{\partial v}{\partial r}\nabla v
    - \frac{1}{2} |\nabla v|^2 \nabla r \right)
    - \frac{1}{2} |\nabla v|^2 & \Delta r
     + (\nabla dr)(\nabla v,\nabla v)  \nonumber \\
    & + 2(\rho '(r) + c)\left( \frac{\partial v}{\partial r} \right)^2
      - qv \frac{\partial v}{\partial r} = 0. 
\end{align}
Since for any vector field $X$
$$
  r^{\gamma}{\rm div}\, X 
  = {\rm div}\,(r^{\gamma}X) - \gamma \,r^{\gamma -1}(Xr),
$$
multiplying the equation $(12)$ by $r^{\gamma }$ further yields
\begin{align*}
  & - \mathrm{div}\left(
    r^{\gamma}\frac{\partial v}{\partial r}
    \nabla v - \frac{1}{2} r^{\gamma} |\nabla v|^2 \nabla r 
    \right)
    + \gamma r^{\gamma-1} \left( \frac{\partial v}{\partial r} \right)^2 
     -\frac{1}{2} r^{\gamma-1}(\gamma +r\Delta r)|\nabla v|^2 \\
  & \hspace{25mm} + r^{\gamma}(\nabla dr)(\nabla v,\nabla v) 
    + 2 ( \rho'(r) + c ) r^{\gamma} 
      \left( \frac{\partial v}{\partial r} \right)^2
   - r^{\gamma} qv \frac{\partial v}{\partial r} = 0.
\end{align*}
Integrating this inequality over $B(s,t)$ with respect to the measure $d\mu_c$, we get  
\begin{align*}
    &  -\left( \int_{S(t)} - \int_{S(s)} \right)
        \left\{ r^{\gamma} \left( \frac{\partial v}{\partial r} \right)^2
       - \frac{1}{2} r^{\gamma} |\nabla v|^2 \right\}\,dA_c \\
    &  + \int_{B(s,t)} r^{\gamma} (\nabla dr)(\nabla v,\nabla v) \,d\mu_c
       + \int_{B(s,t)} \gamma r^{\gamma-1}
       \left( \frac{\partial v}{\partial r} \right)^2 \,d\mu_c\\
    &  - \frac{1}{2} \int_{B(s,t)} r^{\gamma-1} (\gamma + r\Delta r - 2cr)
        |\nabla v|^2 \,d\mu_c \\
    &  + 2 \int_{B(s,t)}\rho '(r) r^{\gamma}
       \left( \frac{\partial v}{\partial r} \right)^2 \,d\mu_c
       - \int_{B(s,t)} r^{\gamma} qv \frac{\partial v}{\partial r} 
       \,d\mu_c = 0
\end{align*} 
by integration-by-parts-lemma. 
Since
\begin{align*}
     2qv \frac{\partial v}{\partial r} r^{\gamma}
    = {\rm div}(r^{\gamma} q v^2 \nabla r)
       -r^{\gamma-1} qv^2 (\gamma +r\Delta r)
       -r^{\gamma} \frac{\partial q}{\partial r} v^2,
\end{align*}
integrating this equation over $B(s,t)$ with respect to $d\mu_c$, the last term of the equation above turns out to be  
\begin{align*}
    & - \int_{B(s,t)} r^{\gamma} qv \frac{\partial v}{\partial r} \,d\mu_c \\
   =& - \frac{1}{2} \left( \int_{S(t)}-\int_{S(s)} \right) 
        r^{\gamma} qv^2 \,dA_c
      + \frac{1}{2} \int_{B(s,t)} r^{\gamma-1} qv^2 
        ( \gamma + r\Delta r- 2cr ) \,d\mu_c \\
    & + \frac{1}{2} \int_{B(s,t)} r^{\gamma} \frac{\partial q}{\partial r} 
        v^2 \,d\mu_c.
\end{align*}
We get Proposition $3.2$ from these two equations:
\begin{align*}
   & \left( \int_{S(t)} - \int_{S(s)} \right) r^{\gamma}
     \left\{ \left( \frac{\partial v}{\partial r} \right)^2
     -\frac{1}{2} |\nabla v|^2 + \frac{1}{2} qv^2 \right\} \,dA_c \\
  =& \int_{B(s,t)} r^{\gamma -1} \left\{ r(\nabla dr)(\nabla v,\nabla v)
     + \gamma \left( \frac{\partial v}{\partial r} \right)^2 \right\} 
       \,d\mu _c  \\
   & - \frac{1}{2} \int_{B(s,t)} r^{\gamma-1} ( \gamma + r\Delta r - 2cr )
     \left( |\nabla v|^2 - qv^2 \right) \,d\mu_c  \\
   &  + \frac{1}{2} \int_{B(s,t)} r^{\gamma} \frac{\partial q}{\partial r} v^2 
      \,d\mu _c
     + 2 \int_{B(s,t)} \rho' (r) r^{\gamma}
     \left( \frac{\partial v}{\partial r} \right)^2 \,d\mu_c  \\
  =& \int_{B(s,t)} r^{\gamma -1} \left\{ r(\nabla dr)(\nabla v,\nabla v)
     - \frac{1}{2} ( \gamma + r\Delta r - 2cr ) \left( |\nabla v|^2 
     - \left( \frac{\partial v}{\partial r} \right)^2 \right) \right\} 
      \,d\mu _c   \\
   & + \int_{B(s,t)} r^{\gamma -1} 
      \left\{ \frac{1}{2} ( \gamma - r\Delta r + 2cr ) + 2r\rho'(r) \right\}
      \left( \frac{\partial v}{\partial r} \right)^2 \,d\mu _c   \\
   &  + \frac{1}{2} \int_{B(s,t)}r^{\gamma -1} 
       \left\{ 
       ( \gamma + r\Delta r - 2cr ) q + r \frac{\partial q}{\partial r} 
     \right\} v^2 \,d\mu _c .
\end{align*}
\end{proof}
From Proposition $3.1$ and $3.2$, we get the following:
%
%
\begin{prop}
Let $\nabla r,X_1,X_2,\cdots ,X_{n-1}$ be an orthonormal base for the tangent space $T_xM$ at each point $x\in M-\overline{U}$. 
Then, for any real numbers $\gamma$, $\varepsilon$, and $0\le s<t$, we have
\begin{align*}
  & \int_{S(t)} r^{\gamma} 
    \left\{ \left( \frac{\partial v}{\partial r} \right)^2
    + \frac{1}{2} q|v|^2 - \frac{1}{2} |\nabla v|^2 
    + \frac{\gamma-\varepsilon }{2r}
      \frac{ \partial v}{\partial r}v \right\} \,dA_c  \\
  & + \int_{S(s)} r^{\gamma} \left\{ \frac{1}{2}|\nabla v|^2
    - \frac{1}{2} q |v|^2 - \left( \frac{\partial v}{\partial r} \right)^2
    - \frac{ \gamma - \varepsilon }{2r}
    \frac{\partial v}{\partial r} v \right\} \,dA_c  \\
 =& \int_{B(s,t)} r^{\gamma -1} \left\{ r(\nabla dr)(\nabla v,\nabla v)
    - \frac{1}{2} ( r \Delta r - 2cr + \varepsilon )
     \sum_{i=1}^{n-1} \left( dv(X_i) \right)^2 \right\} \,d\mu _c\\
  & + \int_{B(s,t)}r^{\gamma-1} 
      \left\{ \gamma - \frac{1}{2} ( r \Delta r - 2cr + \varepsilon ) 
      + 2r \rho'(r) \right\} \left( \frac{\partial v}{\partial r} \right)^2 
      \,d\mu _c \\
  & + \frac{1}{2}\int_{B(s,t)} r^{\gamma-1}
     \left\{ r \left( \frac{\partial q}{\partial r} \right)
     +q(r\Delta r - 2cr + \varepsilon ) \right\} |v|^2 \,d\mu _c \\
  & + \frac{\gamma -\varepsilon}{2} \int_{B(s,t)} r^{\gamma-1}
     \left\{ \frac{\gamma-1}{r} + 2\rho '(r) \right\} 
     \frac{\partial v}{\partial r}v \,d\mu _c.
\end{align*}
\end{prop}
\begin{proof}
When we set $\psi =r^{\gamma -1}$ in Proposition $3.1$, we get 
\begin{align*}
   & \left( \int_{S(t)} - \int_{S(s)} \right) r^{\gamma-1} 
     \frac{\partial v}{\partial r}v \,dA_c\\
 = & \int_{B(s,t)} r^{\gamma -1} \left\{ |\nabla v|^2 - qv^2 \right\} \,d\mu _c
     + \int_{B(s,t)} 
     \Bigl( (\gamma -1) r^{\gamma-2} + 2r^{\gamma-1} \rho'(r) \Bigr)
     \frac{\partial v}{\partial r} v \,d\mu _c.
\end{align*}
If we multiply both sides of this equation by $\frac{\gamma-\varepsilon}{2}$ and add it to the equation in Proposition $3.2$, we get Proposition $3.3$.
\end{proof}
\begin{lem}
We have for any real number $\beta$
\begin{align*}
   & \left( \int_{S(t)} - \int_{S(s)} \right) r^{\beta} v^2 \,dA_c
     \nonumber\\
 = & \int_{B(s,t)} r^{\beta}
     \left\{ \left((\Delta r-2c) + \frac{\beta}{r}\right) v^2
   + 2v \frac{\partial v}{\partial r} \right\} \,d\mu_c.
\end{align*}
\end{lem}
\begin{proof}
A direct computation shows that
\begin{align*}
     \mathrm{div}( r^{\beta}v^2\nabla r)
    = r^{\beta}
    \left\{ \left( \Delta r + \frac{\beta}{r} \right)v^2
    +2v\frac{\partial v}{\partial r}\right\}.
\end{align*}
Integrating this equation with respect to $d\mu_c$ and using the integration-by-parts-lemma, we get Lemma $3.2$. 
\end{proof}
\section{Faster than polynomial decay}

In the following, we shall use the following convention for the sake of simplicity: for any real number $a$, we define $\widehat{a}$ by 
\begin{align*}
   \widehat{a} := (n-1) a.
\end{align*}
The proof of Theorem $1.1$ will be accomplished by following three procedures: (1)to show faster than polynomial decay; (2) to show faster than exponential decay; (3) to show vanishing on a neighborhood of infinity. 
Section $3$, $4$, and $5$ will be devoted to these procedures $(1)$, $(2)$, and $(3)$, respectively. 
\begin{thm}
Let $(M,g)$ be an $n$-dimensional noncompact complete Riemannian manifold and $U$ an open subset of $M$. 
Assume that the end $E := M-U$ has radial coordinates and write $r={\rm dist}\,(U,*)$. 
We also assume that there exists a positive constant $r_0$ such that 
\begin{align}
   \left( 1 - \frac{A_1}{r} \right) \widetilde{g}
  \le 
   \nabla dr \le \left( 1 + \frac{B_1}{r} \right) \widetilde{g}
   \qquad {\rm on}~~B(r_0,\infty)  \tag{$*_1$}
\end{align}
and 
\begin{align}
   {\rm Ric}\,(\nabla r, \nabla r) 
  \ge - (n-1) \left( 1 + \frac{b_1}{r} \right) 
   \qquad {\rm on}~~B(r_0,\infty)  \tag{$*_2$}
\end{align}
where we set $\widetilde{g}= g - dr\otimes dr$, and $A_1$, $B_1$, and $b_1$ are positive constants. 
Let $\lambda >0$ be a constant and $u$ a solution to the following$\,:$
\begin{align}
   L u + \lambda u = 0  \quad\qquad {\rm on}~~B(r_0,\infty) \tag{$8$}
\end{align}
and assume that a constant $\gamma >0$ satisfies
\begin{align}
   \lambda \Big\{ 2\gamma - ( \widehat{A}_1 + \widehat{B}_1) \Big\}
   >
   (n-1) \left( \widehat{A}_1 + \frac{ \widehat{b}_1 }{2} \right). 
   \tag{$*_{3}$}
\end{align}
If $u$ satisfies the condition
\begin{equation}
   \liminf_{t\to \infty}~t^{\gamma}\int_{S(t)}
   \left\{ \left( \frac{\partial u}{\partial r} \right) ^2 
   + |u|^2 \right\} \,dA_c=0,
\end{equation}
then we have for any $m>0$
\begin{equation}
   \int_{B(r_0,\infty)} r^m \left\{ |\nabla u|^2 + |u|^2 \right\} \, d\mu_c
   < \infty.
\end{equation}
\end{thm}
\begin{proof}
First, note that our assumptions $(*_1)$ and $(*_2)$ imply that
\begin{align}
   & r(\nabla dr)(\nabla v,\nabla v)
     \ge (r - A_1) \sum_{i=1}^{n-1} \left( dv(X_i) \right)^2 ; \tag{$*_{4}$} \\
   & -\frac{ \widehat{A}_1 }{r} \overset{(*_{5.1})}{\le} 
     \Delta r - 2c \overset{(*_{5.2})}{\le} \frac{\widehat{B}_1}{r}, 
     \tag{$*_{5}$}
\end{align}
and
\begin{align}
    - \frac{\partial (\Delta r)}{\partial r}
   = & |\nabla dr|^2+\mathrm{Ric}\,(\nabla r,\nabla r) \nonumber \\
 \ge & (n-1)\left(1-\frac{A_1}{r}\right)^2
       -(n-1)\left(1+\frac{b_1}{r}\right) \nonumber \\
   = & -\frac{1}{r} 
       \left( 2\widehat{A}_1 + \widehat{b}_1 - \frac{\widehat{A}_1A_1}{r} \right), \tag{$*_{6}$}
\end{align} 
where the first identity is due to Proposition $2.3$. 
From the assumption $(*_3)$, we can choose a constant $\varepsilon$ so that 
\begin{align}
    2 \gamma - \widehat{B}_1 > \varepsilon > \widehat{A}_1 
    + \frac{(n-1) \left( \widehat{A}_1 + \frac{\widehat{b}_1}{2} \right)}
      {\lambda}.
\end{align}
We shall put $\rho(r)=0$ in Proposition $3.3$; 
then, 
\begin{align*}
   v = u;
\end{align*}
moreover, from $(*_5)$ and $(*_6)$,
\begin{align}
  &  \lambda - \frac{(n-1) \widehat{B}_1}{2r} 
     \le  q = \lambda + c( 2c - \Delta r) 
     \le  \lambda + \frac{(n-1) \widehat{A}_1}{2r};\\
  &  r \, \frac{\partial q}{\partial r}
     = -cr \frac{\partial (\Delta r)}{\partial r} \ge 
     - \frac{(n-1)}{2} 
     \left( 2\widehat{A}_1 + \widehat{b}_1 - \frac{\widehat{A}_1A_1}{r} \right)
     ;\nonumber
\end{align}
and hence, 
\begin{align}
     & r \, \frac{\partial q}{\partial r} 
       + q( r\Delta r -2cr + \varepsilon ) \nonumber \\
 \ge & - \frac{(n-1)}{2} 
       \left( 
       2 \widehat{A}_1 + \widehat{b}_1 - \frac{\widehat{A}_1A_1}{r} 
       \right)
        + \left( \lambda - \frac{(n-1) \widehat{B}_1}{2r} \right)
       \left( \varepsilon - \widehat{A}_1 \right) \nonumber \\
   = & \lambda \varepsilon - \lambda \widehat{A}_1
       - (n-1) \widehat{A}_1 - \frac{(n-1)}{2}\widehat{b}_1 
       + O(r^{-1}) \nonumber \\
   = & C_1 + O(r^{-1}),
\end{align}
where we set $C_1 =  \lambda \varepsilon - \lambda \widehat{A}_1 - (n-1) \widehat{A}_1 - \frac{(n-1)}{2}\widehat{b}_1$. 
Note that $C_1 >0 $ by $(15)$. 

We also have
\begin{align}
   - \left( \frac{\partial u}{\partial r} \right)^2
   - \frac{\gamma - \varepsilon}{2r}\frac{\partial u}{\partial r}u
  \le  \frac{(\gamma - \varepsilon)^2}{16r^2} |u|^2 ,
\end{align}
and by Lemma $3.2$ with $\beta = \gamma -2$ we see that
\begin{align}
   &  \frac{(\gamma - \varepsilon)^2}{16} 
      \left( \int_{S(t)} - \int_{S(s)} \right) r^{\gamma} \frac{|u|^2}{r^2} 
      \, dA_c \nonumber \\
 = &  \int_{B(s,t)} r^{\gamma -1}
      \left\{ 
      O(r^{-2})|u|^2 + O(r^{-1}) u\frac{\partial u}{\partial r} 
      \right\} \,d\mu_c
\end{align}

Substituting $(16)$, $(17)$, $(18)$, and $(19)$ in Proposition $3.3$ with $\rho(r) = 0$, we see that  
\begin{align}
  & \int_{S(t)} r^{\gamma}
    \left\{ \left( \frac{\partial u}{\partial r} \right)^2
    + \frac{1}{2}\bigl( \lambda + O(r^{-1}) \bigr) |u|^2 \right\} 
    \,dA_c \nonumber \\
  & + \int_{S(s)}r^{\gamma}
    \left\{ \frac{1}{2}|\nabla u|^2 - \frac{1}{2}q|u|^2 \right\} 
    \,dA_c \nonumber \\
 \ge
  &  \int_{B(s,t)} r^{\gamma-1} 
     \left\{ r - A_1 - \frac{1}{2}\Bigl( (n-1)B_1 + \varepsilon \Bigr) \right\}
     \sum_{i=1}^{n-1} \left( du(X_i) \right)^2 \,d\mu_c \nonumber \\
  &  + \frac{1}{2} \int_{B(s,t)} r^{\gamma-1}
     \bigl\{ C_2 + O(r^{-1}) \bigr\}
     \left(\frac{\partial u}{\partial r}\right)^2 \,d\mu_c \nonumber \\
  & + \frac{1}{2}\int_{B(s,t)}r^{\gamma-1}
     \left\{ C_1 + O(r^{-1}) \right\}|u|^2 \,d\mu_c,
\end{align}
where we set $C_2 = 2\gamma - \varepsilon- (n-1)B_1$. 
Note that $C_2 > 0$ by $(15)$. 
Therefore, if we take sufficiently large constant $r_1 > 0$, then for any $t > s\ge r_1$ the right hand side of $(20)$ is bounded from below by 
$$
   \frac{1}{4}\min\{C_1,C_2 \} \int_{B(s,t)}r^{\gamma-1}
   \left\{|\nabla u|^2+|u|^2\right\} \,d\mu_c.
$$
Our assumption $(13)$ implies that there exits a divergent sequence $\{t_i\}$ of numbers such that the first term with $t=t_i$ of the inequality above converges to zero as $i\to \infty$. 
Hence, putting $t=t_i$ and letting $i\to \infty$, we get
\begin{align}
   &  \int_{S(s)}r^{\gamma}
     \left\{ |\nabla u|^2-q|u|^2 \right\} \,dA_c \nonumber \\
     \ge
   & \frac{1}{2} \min\{C_1,C_2 \} \int_{B(s,\infty)}
     r^{\gamma-1} \left\{ |\nabla u|^2+|u|^2 \right\} \,d\mu_c
\end{align}
for $s\ge r_1$. 
Integrating this inequality with respect to $s$ over $[t,t_1]$ $(r_1\le t<t_1)$, we have
\begin{align*}
     &  \frac{1}{2} \min\{C_1,C_2 \} 
        \int^{t_1}_t \,ds \int_{B(s,\infty)} r^{\gamma-1}
        \left\{ |\nabla u|^2 + |u|^2 \right\} \,d\mu_c \\
 \le &  \int_{B(t,t_1)} r^{\gamma}
        \left\{ |\nabla u|^2 - q|u|^2 \right\} \,d\mu_c\\
  =  &  \left( \int_{S(t_1)} - \int_{S(t)} \right)
        r^{\gamma} \frac{\partial u}{\partial r}u \,dA_c
        - \gamma \int_{B(t,t_1)} r^{\gamma-1}
        \frac{\partial u}{\partial r}u \,d\mu_c.
\end{align*}
In the last line, we have used the equation in Proposition $3.1$ with $\rho (r)=0$ and $\psi=r^{\gamma}$. 
Since our assumption $(13)$ implies 
$$
    \liminf_{t_1\to \infty}\int_{S(t_1)} r^{\gamma}
    \frac{\partial u}{\partial r}u \,dA_c=0,
$$
letting $t_1\to \infty$ and using Fubini's theorem, we have 
\begin{align}
   &  \min\{C_1,C_2 \} 
      \int^{\infty}_t \,ds \int_{B(s,\infty)}r^{\gamma-1}
      \left\{ |\nabla u|^2 + |u|^2 \right\} \,d\mu_c  \nonumber \\
 = &  \min\{C_1,C_2 \} \int_{B(t,\infty)} (r-t) r^{\gamma-1}
      \left\{ |\nabla u|^2 + |u|^2 \right\} \,d\mu_c \nonumber \\
      \le 
   &  \int_{S(t)} r^{\gamma}
      \left\{ \left( \frac{\partial u}{\partial r} \right)^2 + |u|^2 
      \right\} \,dA_c
      + \gamma \int_{B(t,\infty)} r^{\gamma-1}
      \left\{ \left( \frac{\partial u}{\partial r} \right)^2 + |u|^2 \right\}
      \,d\mu_c < \infty,
\end{align}
where the right hand side of this inequality is finite by $(21)$. 
Hence we see that the desired assertion $(14)$ holds for $m = \gamma$. 

Integrating this inequality $(22)$ with respect to $t$ over 
$[t_1,\infty)~(t_1\ge r_1)$ and using Fubini's theorem, we get
\begin{align*}
   & \min\{C_1,C_2 \} 
     \int_{B(t,\infty)} (r-t)^2 r^{\gamma-1} 
     \{ |\nabla u|^2 + |u|^ 2\} \,d\mu_c \\
     \le 
   & \int_{B(t,\infty)} r^{\gamma}
     \left\{ \left( \frac{\partial u}{\partial r} \right)^2 + |u|^2 \right\} 
     \,d\mu_c
     + \gamma \int_{B(t,\infty)} (r-t) r^{\gamma-1}
     \left\{ \left( \frac{\partial u}{\partial r} \right) ^2 + |u|^2 \right\} 
     \,d\mu_c\\
 < & \infty,
\end{align*}
where the right hand side of this inequality is finite by $(22)$. 
Thus, we see that the desired assertion $(14)$ holds for $m=\gamma+1$. 
Repeating the integration with respect to $t$ shows that the assertion $(14)$ is valid for $m = \gamma + 2, \gamma + 3, \cdots $, therefore, for any $m>0$.   
\end{proof}
\section{Faster than exponential decay}
\begin{lem}
We have for any $\alpha>0$ and $\gamma \in {\rm R}$
\begin{align*}
   & \alpha \left( \int_{S(t)} - \int_{S(s)}\right) r^{\gamma -1} 
     |v|^2 \,dA_c \\
     \ge 
   & - \alpha \int_{B(s,t)} r^{\gamma -1} 
     \left\{ \left( \frac{\partial v}{\partial r} \right)^2 
     + \left( 1 - \frac{ \gamma - 1 - \widehat{A}_1 }{r} \right) |v|^2 
     \right\} \,d\mu_c.
\end{align*}
\end{lem}
\begin{proof}
Set $ \beta = \gamma -1 $ in Lemma $3.2$. 
Then by $(*_{5.1})$, we have 
\begin{align*}
   & \left(\int_{S(t)}-\int_{S(s)}\right) r^{\gamma -1} |v|^2 \,dA_c \\
 = & \int_{B(s,t)} r^{\gamma -1} \left\{ \left( 
     \Delta r - 2c + \frac{\gamma -1}{r} \right) |v|^2
     + 2v \frac{\partial v}{\partial r} 
     \right\} \,d\mu_c \\
     \ge 
   & \int_{B(s,t)} r^{\gamma -1} 
     \left\{ 
     \frac{ \gamma -1 - \widehat{A}_1 }{r} |v|^2
     - |v|^2 - \left( \frac{\partial v}{\partial r} \right)^2 
     \right\} \,d\mu_c \\
 = & - \int_{B(s,t)} r^{\gamma -1} 
     \left\{ \left(\frac{\partial v}{\partial r}\right)^2 
     + \left( 1 - \frac{ \gamma - 1 - \widehat{A}_1 }{r}  \right) |v|^2 
     \right\} \,d\mu_c.
\end{align*}
\end{proof}
\begin{thm}
Under the assumptions of Theorem $4.1$, we have for any $ k > 0 $
\begin{align}
   \int_{B(r_0,\infty)}
   e^{kr} \left\{ u^2 + |\nabla u|^2 \right\} \,d\mu_c < \infty.
\end{align}
\end{thm}
\begin{proof}
In order to prove this theorem, we shall assume that 
\begin{equation}
   \int_{B(r_0,\infty)} r^{m} 
   \left\{ u^2 + |\nabla u|^2 \right\} \,d\mu_c < \infty
   \qquad \mathrm{for~all}~~m \ge 1
\end{equation}
and show that
\begin{align*}
  \int_{B(r_0,\infty)} e^{ k r} 
  \left\{ u^2+|\nabla u|^2 \right\} \,d\mu_c < \infty 
  \qquad {\rm for~any}~~ k > 0.
\end{align*}
For that purpose, let us set 
\begin{align}
   \rho(r) = m \log r
\end{align} 
and $\gamma = \varepsilon$ in Proposition $3.3$. 
Here, $\gamma > 0$ is a large constant determined later. 
Then we have
\begin{align}
  v & = r^m u;\\
  q & =  \lambda + \left( \rho '(r) \right)^2 - \rho ''(r) 
         + (2c - \Delta r)\left( \rho '(r) + c \right) \nonumber \\
    & = \lambda +\left( \frac{m}{r} \right)^2 + \frac{m}{r^2}
        + (2c - \Delta r) \left( \frac{m}{r} + c \right)  \\
    & \ge 
        \lambda -c \frac{\widehat{B}_1}{r} 
        + \frac{m^2}{r^2}\left\{ 1 + \frac{1-\widehat{B}_1}{m} \right\}
        \qquad {\rm by}~~(*_{5.2}); \nonumber \\
      r \frac{\partial q}{\partial r}
    & = - 2\frac{m^2+m}{r^2}
        - r \frac{ \partial (\Delta r) }{ \partial r }
        \left( c + \frac{m}{r} \right)
        + ( \Delta r - 2c ) \frac{m}{r}  \nonumber \\
    & \ge 
        - 2\frac{m^2+m}{r^2} - 
        \left( 2\widehat{A}_1 + \widehat{b}_1 
        - \frac{\widehat{A}_1 A_1}{r} \right) 
        \left( c + \frac{m}{r} \right) 
        - \frac{\widehat{A}_1m}{r^2}  
        \qquad {\rm by}~(*_{5.1})~{\rm and}~(*_6)  \nonumber \\
    & = - c ( 2\widehat{A}_1 + \widehat{b}_1) - 
        \frac{m}{r} 
        \left\{ 
        2\widehat{A}_1 + \widehat{b}_1 - \frac{c\widehat{A}_1A_1}{m} 
        \right\} 
        - \frac{2m^2}{r^2} 
        \left\{ 1 + \frac{ 2 + \widehat{A}_1 (1 - A_1) }{2m} 
        \right\}, \nonumber 
\end{align}
and hence,  
\begin{align*}
  &  r \frac{\partial q}{\partial r} 
     + q \left( r \Delta r -2cr + \gamma \right)\\
 \ge 
  &  - c ( 2 \widehat{A}_1 + \widehat{b}_1 ) 
     - \frac{m}{r} \left\{ 2\widehat{A}_1 + \widehat{b}_1 
     - \frac{c\widehat{A}_1A_1}{m} \right\} 
     - \frac{2m^2}{r^2} \left\{ 1 + \frac{2+\widehat{A}_1(1-A_1) }{2m} 
       \right\}\\
  &  + \left\{ \lambda -  \frac{ c \widehat{B}_1}{r}  
     + \frac{m^2}{r^2}\left( 1 + \frac{1-\widehat{B}_1}{m} \right) \right\}
     \left( \gamma - \widehat{A}_1 \right)\\
 = 
  &  \left( \gamma - \widehat{A}_1 \right) \lambda 
     - c \left( 2\widehat{A}_1 + \widehat{b}_1 \right) 
     - \frac{m}{r} \left\{ 2 \widehat{A}_1 + \widehat{b}_1 
     + \frac{c}{m} 
      \left( \widehat{B}_1 \gamma - \widehat{A}_1( \widehat{B}_1 + A_1 ) \right)     \right\} \\
  &  + \frac{m^2}{r^2} 
     \left\{ 
     \gamma -2 - \widehat{A}_1  
     + \frac{ \left(1 - \widehat{B}_1 \right)
     \left( \gamma - \widehat{A}_1 \right) - 2 - \widehat{A}_1 (1- A_1) }{m} 
     \right\}.
\end{align*}
Thus, for
\begin{align}
   \gamma > \gamma_1( \lambda, A_1, b_1 ) := 
   \max \left\{ \widehat{A}_1 
   + \frac{ c( 2\widehat{A}_1 + \widehat{b}_1  ) + 1 }{\lambda}, 
   2 + \widehat{A}_1 + \frac{\widehat{A}_1 A_1}{\widehat{B}_1} \right\}  
\end{align}
and $ 0 < \alpha < \frac{1}{2} $, we set 
\begin{align*}
       \widehat{\lambda}_1 
  := & \left( \gamma - \widehat{A}_1 \right) \lambda 
       - c \left( 2\widehat{A}_1 + \widehat{b}_1 \right)
       \hspace{20mm} \bigl(\, \ge 1~~{\rm by}~(28) \bigr);\\
       \widehat{\lambda}_2 
  := & \widehat{\lambda}_1 
       - 2 \alpha  + 2 \alpha \frac{ \gamma - 1 - \widehat{A}_1 }{r}  
       \hspace{25.5mm} 
       \bigl(\, \ge  1 - 2 \alpha > 0 ~~ {\rm by}~(28) \bigr) ; \\
       \widehat{C}_1
  := & 2 \widehat{A}_1 + \widehat{b}_1 
       + \frac{c}{m}
       \left( \widehat{B}_1 \gamma - 
       \widehat{A}_1( A_1 + \widehat{B}_1 ) \right)
       \hspace{6.5mm} \bigl(\, >0~~{\rm by}~(28) \bigr); \\
      \widehat{C}_2
  := & \gamma -2 - \widehat{A}_1  
       + \frac{ \left(1 - \widehat{B}_1 \right) 
       \left( \gamma - \widehat{A}_1 \right) - 2 
       - \widehat{A}_1 (1- A_1) }{m} \\
   = & \gamma \left( 1 + \frac{1 - \widehat{B}_1}{m} \right) -2 - \widehat{A}_1
       - \frac{\left(1 - \widehat{B}_1 \right)\widehat{A}_1 
       + 2 + \widehat{A}_1 (1- A_1) }{m}.
\end{align*}
What is more, by taking $m_1(B_1) \ge 1$ and $\gamma_2( \lambda, A_1, B_1, b_1) \ge \gamma_1( \lambda, A_1, b_1 )$ sufficiently large, we see that 
\begin{align}
    \widehat{C}_2 ~\left( \ge \frac{1}{2} \gamma \right) > 0 \qquad 
    {\rm for}~~m \ge m_1(B_1)~
    {\rm and}~\gamma \ge \gamma_2 ( \lambda, A_1, B_1, b_1).
\end{align}
Thus, we get 
\begin{align*}
     & r \frac{\partial q}{\partial r} 
       + q \left( r \Delta r -2cr + 1 \right)
       - 2 \alpha \left( 1 - \frac{\gamma - 1 - \widehat{A}_1 }{r} \right)\\
 \ge & \widehat{\lambda}_2 
       - \frac{m}{r} \widehat{C}_1 + \frac{m^2}{r^2} \widehat{C}_2.
\end{align*}
By $(*_{5.2})$, we also have
\begin{align*}
   \gamma - \frac{1}{2}( r \Delta r - 2cr + \gamma ) + 2r\rho'(r)
   \ge 
   2m + \frac{1}{2} \left( \gamma - \widehat{B}_1 \right).
\end{align*}
Therefore, combining Lemma $5.1$, Proposition $3.3$ and estimates above yield
\begin{align}
  & \int_{S(t)} r^{\gamma} 
    \left\{ \left( \frac{\partial v}{\partial r} \right)^2
    + \frac{1}{2} q|v|^2 + \frac{\alpha }{r} |v|^2 \right\} \,dA_c\nonumber \\
  & + \frac{1}{2} \int_{S(s)} r^{\gamma} 
    \left\{|\nabla v|^2 - q|v|^2 \right\} \,dA_c
    - \int_{S(s)} r^{\gamma} 
    \left\{ \left( \frac{\partial v}{\partial r} \right)^2
    + \frac{\alpha }{r}|v|^2 \right\} \,dA_c\nonumber \\
    \ge 
  & \int_{B(s,t)} r^{\gamma -1} \left\{ \left(r - \widehat{A}_1 \right) - 
     \frac{1}{2} \left( \gamma + \widehat{B}_1 \right) \right\}
    \sum_{i=1}^{n-1} \left( dv(X_i) \right)^2 \,d\mu_c\nonumber \\
  & + \int_{B(s,t)} r^{\gamma -1} 
    \left\{ 2m + \frac{1}{2} \left( \gamma - \widehat{B}_1 \right) 
    - \alpha  \right\}
    \left( \frac{\partial v}{\partial r} \right)^2 \,d\mu_c\nonumber \\
  & + \frac{1}{2} \int_{B(s,t)} r^{\gamma -1} 
    \left\{ \widehat{\lambda}_2 
    - \widehat{C}_1\frac{m}{r} + \widehat{C}_2\frac{m^2}{r^2} \right\}
    |v|^2 \,d\mu_c .
\end{align}
In view of the right hand side of this inequality, we see that the first term is nonnegative for any $t > s \ge r_1( A_1, B_1 , \gamma ) := \frac{1}{2}(\gamma + \widehat{B}_1 ) + \widehat{A}_1 > r_0 $; also, the second term is nonnegative if $m \ge \frac{1+\widehat{B}_1}{4}$. 
Besides, according $(27)$, we see that $\lim_{r\to \infty} q = \lambda$, and hence, $(26)$ and our assumption $(24)$ imply that 
$$
   \liminf_{t\to \infty} \int_{S(t)} r^{\gamma} 
   \left\{ \left( \frac{\partial v}{\partial r} \right)^2 
   + \frac{1}{2} q|v|^2 + \frac{\alpha }{r}|v|^2 \right\} \,dA_c=0.
$$
Therefore, taking an appropriate divergent sequence $\{t_i\}$ of numbers, putting $t=t_i$ in $(30)$ and letting $i\to \infty$, we get, for any $m \ge m_2(B_1):= \max \left\{ \widehat{B}_1 + 1, m_1(B_1) \right\}$, $s \ge r_1( A_1,B_1 , \gamma )$, $\gamma \ge \gamma_2( \lambda, A_1, B_1, b_1 )$, and $ 0 < \alpha < \frac{1}{2} $,
\begin{align*}
  & \int_{S(s)} r^{\gamma} \left\{|\nabla v|^2 - q|v|^2 \right\} \,dA_c
    - 2 \int_{S(s)} r^{\gamma} 
    \left\{\left(\frac{\partial v}{\partial r}\right)^2
    + \frac{\alpha }{r}|v|^2\right\} \,dA_c\\
    \ge 
  & \int_{B(s,\infty)} r^{\gamma -1}
    \left\{ \widehat{\lambda}_2 
    - \widehat{C}_1\frac{m}{r} + \widehat{C}_2\frac{m^2}{r^2} \right\}
    |v|^2 \,d\mu_c.
\end{align*}

Multiplying the both sides of this inequality by $s^{1 - 2m -\gamma }$ and integrating it with respect to $s$ over $[x,\infty)$ $\left( x \ge r_1( A_1, B_1 , \gamma) \right)$, we get 
\begin{align}
    &  \int_{B(x,\infty)} r^{1-2m} 
      \bigl\{ |\nabla v|^2 - q|v|^2 \bigr\} \,d\mu_c
      - 2 \int_{B(x,\infty)} r^{1-2m}
      \left\{ \left(\frac{\partial v}{\partial r}\right)^2
      +\frac{\alpha }{r}|v|^2 \right\} \,d\mu_c \nonumber \\
\ge & \int_x^{\infty} s^{1 - 2m -\gamma } \,ds 
     \int_{B(s,\infty)} r^{\gamma -1} 
     \left\{ 
     \widehat{\lambda}_2 
     - \widehat{C}_1\frac{m}{r} + \widehat{C}_2\frac{m^2}{r^2} 
     \right\} 
     |v|^2 \,d\mu_c.
\end{align}
On the other hand, Proposition $3.1$ with $\psi = r^{1-2m}$ implies 
\begin{align}
   & \int_{B(x,\infty)} r^{1-2m} \bigl\{ |\nabla v|^2 - q|v|^2 \bigr\} \,d\mu_c
     \nonumber \\
 = & - \int_{S(x)} r^{1-2m} \frac{\partial v}{\partial r}v \,dA_c
     - \int_{B(x,\infty)} r^{-2m} \frac{\partial v}{\partial r} v \,d\mu_c,
\end{align}
where we have used $(24)$ and $(26)$. 
If we write $dv_g = \sqrt{G}(r,y) \,dr dv_{\partial U}$ ($y \in \partial U$ and $dv_{\partial U}$ is the induced measure on $\partial U$ from the Riemannian measure $dv_g$), then $\frac{\partial \left( \sqrt{G} \right)}{\partial r} = (\Delta r) \,\sqrt{G}$ (this identity follows from the definition of the Riemannian measure and Laplacian; see, for example, \cite{A-K} pp.7.), and hence, a direct computation shows that 
\begin{align}
     -\int_{S(x)} r^{1-2m} \frac{\partial v}{\partial r}v \,dA_c
 = & - \frac{1}{2} \frac{d}{dx} \left( x^{1-2m} \int_{S(x)}|v|^2 \,dA_c \right)
     \nonumber\\
   & - \frac{1}{2}\int_{S(x)} r^{-2m} \bigl\{ 2m-1-r(\Delta r-2c) \bigr\}
     |v|^2 \,dA_c.
\end{align}
From $(32)$, $(33)$, and $(*_{5.2})$, it follows that 
\begin{align*}
    & \int_{B(x,\infty)} r^{1-2m} 
      \left\{ |\nabla v|^2 - q|v|^2 \right\} \,d\mu_c \\
  = & -\frac{1}{2} \frac{d}{dx}\left( x^{1-2m} \int_{S(x)} |v|^2 \,dA_c \right)
      - \frac{1}{2} \int_{S(x)} r^{-2m} 
      \bigl\{ 2m-1-r(\Delta r-2c) \bigr\} |v|^2 \,dA_c \\
    & - \int_{B(x,\infty)} r^{-2m} \frac{\partial v}{\partial r} v \,d\mu_c \\
\le & -\frac{1}{2}\frac{d}{dx}\left( x^{1-2m} \int_{S(x)} |v|^2 \,dA_c \right)
      -\frac{1}{2}\int_{S(x)} r^{-2m} 
      \bigl\{ 2m - 1 - \widehat{B}_1 \bigr\}|v|^2 \,dA_c \\
    & + \int_{B(x,\infty)} r^{-2m} 
      \left\{ 
      \frac{1}{8\alpha} \left( \frac{\partial v}{\partial r} \right)^2 
      + 2 \alpha |v|^2 \right\} \,d\mu_c.
\end{align*}
Substituting this inequality into $(31)$, we get 
\begin{align*}
   & -\frac{1}{2}\frac{d}{dx} \left( x^{1-2m}\int_{S(x)} |v|^2 \,dA_c \right)
     - m \int_{S(x)} r^{-2m}  \left\{ 1- \frac{1 + \widehat{B}_1 }{2m} \right\}
     |v|^2 \,dA_c \\
   & -\int_{B(x,\infty) }r^{1-2m}
     \left\{ 2-\frac{1}{8\alpha r} \right\}
     \left( \frac{\partial v}{\partial r} \right)^2 \,d\mu_c \\
     \ge 
   & \int_x^{\infty}  
     s^{1 - 2m -\gamma } \,ds 
     \int_{B(s,\infty)} r^{\gamma -1} 
     \left\{ 
     \widehat{\lambda}_2 
     - \widehat{C}_1\frac{m}{r} + \widehat{C}_2\frac{m^2}{r^2} \right\} 
     |v|^2 \,d\mu_c.
\end{align*}
Hence, if we set $r_2 := \max \left\{ r_1( A_1, B_1 , \gamma ), \frac{1}{16\alpha} \right\}$, then 
\begin{align}
   & - \frac{1}{2} \frac{d}{dx} \left( x^{1-2m}\int_{S(x)}|v|^2 \,dA_c \right)
     - \frac{m}{2} x^{-2m} \int_{S(x)}|v|^2 \,dA_c  \nonumber \\
     \ge 
   & 
     \int_x^{\infty}  
     s^{1 - 2m -\gamma } \,ds 
     \int_{B(s,\infty)} r^{\gamma -1} 
     \left\{ 
     \widehat{\lambda}_2 - \widehat{C}_1\frac{m}{r} + 
     \widehat{C}_2\frac{m^2}{r^2} 
     \right\} 
     |v|^2 \,d\mu_c
\end{align}
for any $m$, $x$ and $\gamma$ satisfying
\begin{align*}
   m \ge m_2(B_1) ,~
   x >   r_2( A_1, B_1 , \gamma, \alpha),~{\rm and}~
   \gamma > 
   \gamma_2( \lambda, A_1, ,B_1, b_1 ).
\end{align*}

Now, we shall show that the right hand side of $(34)$ is nonnegative for sufficiently large $m$, $x$, and $\gamma$. 
For that purpose, let us consider the following quadratic equation: 
\begin{align}
   \widehat{C}_2 \,y^2 - \widehat{C}_1 \,y + \widehat{\lambda}_2 = 0
\end{align}
and calculate its discriminant $D$. 
Then,
\begin{align}
  D = & \left( \widehat{C}_1 \right)^2 - 4 \widehat{C}_2 \widehat{\lambda}_2 
        \nonumber \\
    = & \left\{ 
        \frac{c \widehat{B}_1}{m} \gamma
        + \left( 2 - \frac{c ( A_1 + \widehat{B}_1)}{m} \right)
        \widehat{A}_1 + \widehat{b}_1 
        \right\} ^2 \nonumber \\
      & - 4 \left\{ \left( 1 + \frac{ 1 - \widehat{B}_1 }{m} \right) \gamma 
        - 2 - \frac{2}{m} 
        - \widehat{A}_1 \left( 1 + \frac{ 1 - \widehat{B}_1 }{m}
        +  \frac{ 1 - A_1 }{m} \right) 
        \right\} \nonumber \\
      & \hspace{7mm} \times 
        \left\{
        \left( \lambda + \frac{2\alpha}{r} \right) \gamma 
        + \left( 2c + \frac{2\alpha}{r} - \lambda \right) \widehat{A}_1 
        - c \widehat{b}_1 - 2 \alpha \left( 1 + \frac{1}{r} \right)
        \right\}.
\end{align}
This equation $(36)$ is a quadratic equation of $\gamma$; the coefficient $H$ of $\gamma ^2$ is calculated as follows: 
\begin{align*}
      H := 
  &   \left( \frac{c \widehat{B}_1}{m} \right)^2 
      - 4 \left( 1 + \frac{ 1 - \widehat{B}_1 }{m} \right)
      \left( \lambda + \frac{2\alpha}{r} \right) \\
  \le 
  &   \left( \frac{c \widehat{B}_1}{m} \right)^2 
      - 4 \left( 1 + \frac{ 1 - \widehat{B}_1 }{m} \right) \lambda 
\end{align*}
Therefore, by taking $m_3( \lambda, B_1 ) \ge m_2(B_1)$ sufficiently large, we see that
\begin{align}
   H < - 3 \lambda \qquad 
   {\rm if}~~m \ge m_3( \lambda, B_1 )~~{\rm and}~~
   r \ge r_0.
\end{align}
Thus, in view of $(36)$ and $(37)$, we see that there exist positive constants 
\begin{align*}
   m_4 = m_4(\lambda, A_1, B_1 )  \ge m_3( \lambda, B_1 ),~~
   r_3 = r_3(\lambda, A_1, B_1 , \gamma , \alpha)  
   \ge r_2(A_1 , B_1 , \gamma , \alpha) ,
\end{align*}
and
\begin{align*}
   \gamma_3 = \gamma_3(\lambda, A_1, B_1, b_1 ) 
   \ge \gamma_2( \lambda , A_1, B_1, b_1). 
\end{align*}
such that
\begin{align*}
    D < 0 \qquad {\rm if}~~m \ge m_4,~r \ge r_3,~{\rm and}~\gamma \ge \gamma_3.
\end{align*}
Hence, by $(29)$, \rm if $m \ge m_4$, $r \ge r_3$, and $\gamma \ge \gamma_3$, then
\begin{align}
   \widehat{C}_2 \,y^2 - \widehat{C}_1 \,y + \widehat{\lambda}_2 > 0 \qquad 
   {\rm for~~any}~~y \in {\bf R}.
\end{align}
Combining $(34)$ and $(38)$, we obtain
\begin{align}
   - \frac{1}{2} \frac{d}{dx} \left( x^{1-2m}\int_{S(x)}|v|^2 \,dA_c \right)
   - \frac{m}{2} x^{-2m} \int_{S(x)}|v|^2 \,dA_c  \ge 0 
\end{align}
for $m \ge m_4( \lambda, A_1 ,B_1 )$, and $x \ge r_3( \lambda, A_1, B_1 , \gamma , \alpha)$. 
Note that the left hand side of this inequality is independent of $\gamma$. \\
Now, let us set 
$$
     F(x) = x^{1-2m} \int_{S(x)} |v|^2 \,dA_c = x \int_{S(x)} |u|^2 \,dA_c.
$$
Then, the left hand side of $(39)$ is equal to 
$$
     -\frac{1}{2}\left( F'(x) + \frac{m}{x}F(x) \right),
$$
and hence, 
\begin{align*}
   F'(x) + \frac{m}{x}F(x) \le 0
\end{align*}
for $m \ge m_4 = m_4( \lambda, A_1 ,B_1 )$, and $x \ge r_3 = r_3(\lambda, A_1, B_1, \alpha)$. \\
For any $k>\frac{m_4}{r_3}$, we set 
\begin{align*}
   m = k x .
\end{align*}
Then, we have
$$
   F'(x) + k F(x) \le 0
$$
for any $x \ge r_3 = r_3(\lambda, A_1, B_1, \alpha)$. 
Set $ G(x) = e^{k x}  F(x)$. 
Then $ G'(x) \le 0 $ for $ x \ge r_3 $, and hence, $G(x) \le G(r_3)$, that is, 
\begin{align*}
    \int_{S(x)} |u|^2 \,dA_c \le x^{-1} e^{-kx} G(r_3) \qquad 
    {\rm for}~~x\ge r_3.
\end{align*}
Thus, we obtain
\begin{equation}
    \int_{B(r_0,\infty)} e^{ kr }|u|^2 \,d\mu _c < \infty \qquad 
    {\rm for~any}~~ k > 0.
\end{equation}

Next, we shall show that $(8)$, $(40)$, and boundedness of $\Delta r$ imply that$$
    \int_{B(r_0,\infty)}e^{ kr }|\nabla u|^2 \,d\mu _c < \infty \qquad 
    {\rm for~any}~~ k > 0.
$$
First, consider the integral
$$
    g(R)
    = 2 \int_{B(r_0,R)} e^{kr} u \frac{\partial u}{\partial r} \,d\mu_c.
$$ 
Then, the integration-by-parts-lemma yields
\begin{align*}
      g(R)
  = & \frac{1}{k} \int_{B(r_0,R)}
      \left\langle 
      \nabla \left( e^{ kr } \right),\nabla \left( u^2 \right) 
      \right\rangle \,d\mu_c \\
  = & \frac{1}{k} \left( \int_{S(R)}-\int_{S(r_0)} \right)
      e^{kr} |u|^2 \,dA_c 
      - \int_{B(r_0,R)} ( \Delta r - 2c + k ) e^{kr}|u|^2 \,d\mu_c.
\end{align*}
Since $\lim_{r\to \infty }(\Delta r-2c) = 0$, $(40)$ implies the existence of the limit, $\lim_{R\to \infty }g(R)$. 
In particular, 
\begin{equation}
   \liminf_{R\to \infty } e^{kR}
   \left| \int_{S(R)} u \frac{\partial u}{\partial r} \,dA_c \right| = 0.
\end{equation}
In Proposition $3.1$, we put $\rho =0$ and $\psi = e^{kr}$, and set $ q_0 = \lambda + c ( 2c - \Delta r)$. 
Then $v=u$, and 
\begin{align*}
    &  \int_{B(r_0,R)} \left\{ |\nabla u|^2 - q_0 |u|^2 \right\} 
       e^{kr} \,d\mu _c\\
 =  &  \left( \int_{S(R)} - \int_{S(r_0)} \right)
       \frac{\partial u}{\partial r} u e^{kr} \,dA_c
       - k \int_{B(r_0,R)}
       e^{kr} \frac{\partial u}{\partial r} u \,d\mu _c\\
\le & \left( \int_{S(R)} - \int_{S(r_0)} \right)
      \frac{\partial u}{\partial r} u e^{kr} \,dA_c
      + \frac{1}{2}\int_{B(r_0,R)}
      e^{kr} \left\{ |\nabla u|^2 + k^2 |u|^2 \right\} \,d\mu _c.
\end{align*}
Hence,
\begin{align*}
     & \frac{1}{2} \int_{B(r_0,R)} e^{kr} |\nabla u|^2 \,d\mu _c \\
 \le & \left(\int_{S(R)}-\int_{S(r_0)}\right)
       \frac{\partial u}{\partial r} u e^{kr} \,dA_c
       + \int_{B(r_0,R)}
       \left\{ \frac{k^2}{2} + q_0 \right\} e^{kr} |u|^2 \,d\mu _c.
\end{align*}
Therefore, $(41)$ and $(42)$ imply that
\begin{equation}
     \int_{B(r_0,\infty)} e^{kr} |\nabla u|^2 \,d\mu _c < \infty
\end{equation}
where $k>0$ is arbitrary. 
Thus, $(40)$ and $(42)$ imply our desired result. 
\end{proof}
\section{Vanishing on a neighborhood of infinity}
%
\begin{thm}
Under the assumption of Theorem $5.1$, there exists a positive constant $r_5$ such that
\begin{align*}
   u\equiv 0\qquad \mathrm{on}~~B(r_5,\infty).
\end{align*}
\end{thm}
\begin{proof}
For any fixed $k\ge 1$ and  
\begin{align}
   \gamma = \varepsilon \ge 2 \widehat{A}_1,
\end{align}
Set $\rho (r) = kr$ in Proposition $3.3$. 
Then, we have
\begin{align}
   v = & e^{kr} u ; \\
   q = & \lambda + |\rho '(r)|^2 - \rho ''(r) 
         + ( 2c - \Delta r )( \rho '(r) + c) \nonumber  \\
     = & \lambda + k^2 +  ( 2c - \Delta r )( k + c )  \\
   \ge & \lambda + k^2 - ( k + c ) \frac{\widehat{B}_1}{r}  
         \hspace{38.5mm} {\rm by}~(*_{5.2}) ; \nonumber \\
         \frac{\partial q}{\partial r}
     = & - ( k + c ) \frac{\partial (\Delta r)}{\partial r} \nonumber \\
   \ge & - \frac{( k + c )}{r} \left( 2 \widehat{A}_1 + \widehat{b}_1 
         - \frac{ \widehat{A}_1 A_1 }{r} \right) \hspace{20mm} {\rm by}~(*_{6})
         \nonumber .
\end{align}
Hence, 
\begin{align}
     & r \frac{\partial q}{\partial r} + q( \gamma + r \Delta r - 2cr )
       \nonumber \\
 \ge & -( k + c ) \left( 2 \widehat{A}_1 + \widehat{b}_1 
       - \frac{ \widehat{A}_1 A_1 }{r} \right)
       + \left( \lambda + k^2 - ( k + c ) \frac{\widehat{B}_1}{r}\right)
       ( \gamma - \widehat{A}_1 ) \nonumber \\
   = & ( \gamma - \widehat{A}_1 ) k^2 
       - \left\{ 2 \widehat{A}_1 + \widehat{b}_1 
       + \frac{ ( \gamma - \widehat{A}_1 )\widehat{B}_1 
       - \widehat{A}_1 A_1 }{r} \right\} k  \nonumber \\
     & \hspace{25mm} + \left( \gamma - \widehat{A}_1 \right)
       \left( \lambda - c \frac{\widehat{B}_1}{r} \right)
       - c \left( 
       2 \widehat{A}_1 + \widehat{b}_1 - \frac{ \widehat{A}_1 A_1 }{r} 
       \right).
\end{align} 
Moreover, we have
\begin{align}
   \frac{\gamma}{2} - \frac{1}{2} ( r \Delta r - 2cr ) + 2kr 
   \ge 2kr + \frac{\gamma}{2} - \frac{ \widehat{B}_1 }{2} 
\end{align}
and 
\begin{align}
   &  r(\nabla dr)(\nabla v,\nabla v)
      - \frac{1}{2} \bigl( r \Delta r - 2cr + \gamma \bigr)
      \sum_{i=1}^{n-1}\left( dv(X_i) \right)^2 \nonumber \\
      \ge  
   &  \left\{ ( r - A_1 ) - \frac{1}{2} ( \widehat{B}_1 + \gamma ) \right\}
      \sum_{i=1}^{n-1} \left( dv(X_i) \right)^2.
\end{align}
In view of $(43)$, we see that there exist constants $k_1 = k_1( A_1, B_1, b_1, ,\gamma, r_0 )$ and $r_5 = r_5 ( A_1, B_1, \gamma )$ such that the right hand sides of $(46)$, $(47)$, and $(48)$ are nonnegative for $ k \ge k_1$ and $ r \ge r_5$. 
Therefore, substituting $(46)$, $(47)$, and $(48)$ in Proposition $3.3$ we have
\begin{align}
    \int_{S(t)} r^{\gamma} 
    \left\{ \left( \frac{\partial v}{\partial r} \right)^2 
    + \frac{1}{2} q|v|^2 \right\} \,dA_c 
    + \int_{S(s)} r^{\gamma} \left\{ \frac{1}{2}|\nabla v|^2
    - \left( \frac{\partial v}{\partial r} \right)^2 \right\} \,dA_c  \ge  0
\end{align}
for $ k \ge k_1 $ and $ t > s \ge r_5 $. 
Note that $(23)$ and $(44)$ imply that
\begin{align*}
    \liminf _{t\to \infty} \int_{S(t)} 
    \bigl\{ |\nabla v|^2 + |v|^2 \bigr\} \,dA_c = 0.
\end{align*}
Also, by $(45)$ and $(*_5)$, 
\begin{align*}
    \lim_{ r \to \infty} q = \lambda + k^2 .
\end{align*}
Hence, taking an appropriate divergent sequence $\{t_i\}$ of numbers, putting $t=t_i$ in $(49)$ and letting $i\to \infty$, we get for any $k\ge k_1$ and $s \ge r_5$
\begin{align*}
     \int_{S(s)} \left\{ \frac{1}{2} |\nabla v|^2
     -\left( \frac{\partial v}{\partial r} \right)^2 \right\} \,dA_c
     \ge 0
\end{align*}
The substitution $v = e^{kr}u $ in this inequality yields 
$$
     e^{2ks} \bigl\{ k^2 I_1(s) + k I_2(s) + I_3(s) \bigr\} \ge 0
$$
for any $k\ge k_1$ and $s\ge r_5$, 
where 
$$
   I_1(s) = - \frac{1}{2} \int_{S(s)} |u|^2 \,dA_c
$$
and $I_2(s)$ and $I_3(s)$ are independent of $k$. 
Therefore $I_1(s)=0$ for $s\ge r_5$. 
That is, $u\equiv 0$ on $B(r_5,\infty )$. 
\end{proof}
Theorem $6.1$ and unique continuation theorem imply that $u\equiv 0$ on $B(r_0,\infty )$.
\section{Proof of Theorem $1.1$ and $1.5$}
We shall prove Theorem $1.1$ by contradiction. 
If there exists $\gamma_0>0$ such that 
$$
  \liminf_{t\to \infty} ~ t^{\gamma_0} \int_{S(t)}
  \left\{ 
  \left(\frac{\partial f}{\partial r}\right)^2 + |f|^2 
  \right\} \,dA < \infty,
$$ 
then for $\gamma_1\in (0,\gamma_0)$ 
\begin{equation}
  \liminf_{t\to \infty} ~ t^{\gamma_1} \int_{S(t)}
  \left\{ 
  \left(\frac{\partial f}{\partial r}\right)^2 + |f|^2 
  \right\} \,dA = 0.
\end{equation}
As is seen in Section $3$, if we set $ u = e^{cr} f$, then $u$ satisfies the equation 
\begin{align}
   Lu + \lambda u = 0 , 
\end{align}
with $ \lambda = \alpha - c^2 = \alpha - \frac{(n-1)^2}{4} > 0$. 
In addition, 
\begin{align*}
  & \left( \frac{\partial f}{\partial r} \right)^2 + |f| ^2
    = \left\{ \left(\frac{\partial u}{\partial r}\right)^2
    - 2cu \frac{\partial u}{\partial r} + ( c^2 + 1 )|u|^2 \right\} e^{-2cr}\\
    \ge 
  & \left\{ \varepsilon \left( \frac{\partial u}{\partial r} \right)^2
    + \left(1-\frac{\varepsilon c^2}{1-\varepsilon } \right)
    |u|^2 \right\} e^{-2cr},
\end{align*}
for any $\varepsilon \in (0,1)$. 
Therefore, if we choose $\varepsilon >0$ sufficiently small, 
the right hand side of this inequality is bounded from below by 
$$
  c(n) \left\{ \left( \frac{\partial u}{\partial r} \right)^2 + |u|^2 
  \right\} e^{-2cr},
$$
where $c(n)>0$ is a constant depending only on $n$. 
Hence by $(50)$ we have 
\begin{equation}
    \liminf_{t\to \infty} ~ t^{\gamma_1} \int_{S(t)}
    \left\{ 
    \left( \frac{\partial u}{\partial r} \right)^2 + |u|^2 
    \right\} \,dA_c = 0.
\end{equation}
In view of $(51)$ and $(52)$, Theorem $6.1$ implies that $u \equiv 0$ on $B(r_0,\infty)$, and hence, $f = e^{-cr} u \equiv 0$ on $ B(r_0,\infty)$. 
Unique continuation theorem implies $f\equiv 0$ on $M-\overline{U}$. 
This contradicts our assumption that $f$ is nontrivial, and hence, Theorem $1.1$ follows.\vspace{2mm}\\
\noindent {\it Proof of Theorem $1.5$} \\
Let $f$ be an $L^2(M,dv_g)$-eigenfunction of $-\Delta $ with eigenvalue $\alpha > \frac{(n-1)^2}{4}$. 
Then $f,|\nabla f|\in L^2(M,dv_g)$ and, in particular, $f,|\nabla f| \in L^2(M-\overline{U},dv_g)$. 
But our growth property in Theorem $1.1$ forces $f$ to be identically zero on $M-\overline{U}$, and by unique continuation theorem it must be zero on $M$. 
Indeed, if 
\begin{equation*}
   t^{\gamma}\int_{S(t)} \left\{ \left( \frac{\partial f}{\partial r} \right) ^2   + |f|^2 \right\} \,dA 
   \ge c_1 > 0 \qquad {\rm on}~~B(r_1,\infty),
\end{equation*}
then 
\begin{equation*}
   \int_{B(r_1,\infty)} 
   \left\{ \left( \frac{\partial f}{\partial r} \right) ^2 
   + |f|^2 \right\} \,dA  \ge \int_{r_1}^{\infty} \frac{1}{t^{\gamma}} \, dt.
\end{equation*}
The right hand side of this inequality is infinity if $\gamma \in (0,1]$. 
Thus $\alpha > \frac{(n-1)^2}{4}$ is not eigenvalue of $-\Delta $ and we have completed the proof of Theorem $1.5$.
\section{Example}
In this section, we shall construct a Riemannian manifold which shows that the curvature decay condition $K + 1 = o(r^{-1})$ mentioned in Theorem $1.4$ and Theorem $1.6$ is sharp. 
In order to do so, we shall use the following theorem essentially due to Atkinson \cite{A}; for the proof of Lemma $8.1$, see Arai-Uchiyama \cite{A-U} and references there. 
\begin{lem}
Let $\lambda > 0 $, $\varepsilon > 0 $, and $k \in {\bf R}$  be constants and 
$q(x) \in C^0[0,\infty)$ a real-valued function. 
Assume that 
\begin{align*}
    q(x) = - k \, \frac{\sin 2x}{x} + O \left( x^{-1-\varepsilon} \right)
    \qquad ( x \to \infty)
\end{align*}
and consider the eigenvalue equation
\begin{align}
    \left( - \frac{d^2}{dx^2} + q(x) \right) w (x) = \lambda w(x) 
    \qquad {\rm on}~~[0,\infty).
\end{align}
Then, the following properties $(a)$ and $(b)$ are eqivalent$\,:$
\begin{enumerate}[$(a)$]
   \item  The equation $(53)$ has a nontrivial solution $w\in L^2[0,\infty)$$;$
   \item  $|k|>2$ and $\lambda = 1 $.
\end{enumerate}
\end{lem}
\vspace{5mm}
\noindent {\bf {\it Proof of Theorem $1.8$ $($Construction of a Riemannian manifold$)$}.}
\vspace{1mm}

Firstly, except for a positive multiplier, we shall construct a desired function $f$ on a neighborhood of infinity as follows: for $r\ge 1$, let us set 
\begin{align}
    f_1 (r)
    = \exp 
    \left\{ \int_{1}^r \left( 1 + k\frac{\sin 2x}{x} \right) \,dx \right\}
\end{align}
and consider the Riemannian manifold with boundary:
\begin{align*}
     (N, g_N) 
     = \Bigl( [1,\infty ) \times S^{n-1}(1) 
     , dr^2 + f_1^{\,2}(r)g_{S^{n-1}(1)} \Bigr),
\end{align*}
where $k \in {\bf R}$ is a constant satisfying
\begin{align}
  |k|(n-1)\sqrt{(n-1)^2 + 4} > 4 . 
\end{align}
Moreover, we set
\begin{align*}
     S(r)   := & \frac{f_1'(r)}{f_1(r)};\\ 
     K(r)   := & -\frac{f_1''(r)}{f_1(r)};\\
     q_0(r) := & \frac{(n-1)(n-3)}{4}S(r)^2 - \frac{(n-1)}{2}K(r).
\end{align*}
Then, direct computations show that 
\begin{align}
    S(r)   = & 1 + \frac{k\sin 2r}{r};\\
    K(r)   = & - 1 -\frac{2\sqrt{2}k\sin (2r+\frac{\pi}{4})}{r}+O(r^{-2});\\
    q_0(r) = & \frac{(n-1)^2}{4 }+ 
               \frac{k(n-1)\sqrt{(n-1)^2+4}}{2r}\sin (2r+c_n)+O(r^{-2}), 
               \nonumber
\end{align}
where $c_n$ is a constant depending only on $n$. 
Hence, in view of $(55)$, Lemma $8.1$ implies that there exists a nontrivial solution $w (x) \in L^2([1,\infty),dx)$ to the equation
\begin{align*}
     \left( -\frac{d^2}{dx^2} + q_0(x) - \frac{(n-1)^2}{4} \right) w (x)
     = w (x).
\end{align*}
Note that $w$ {\it oscillates around} $0$ and $\lim_{x \to \infty} |w|(x) = 0 $ (see \cite{A} and \cite{A-U}). 
Using this function $w$, we define a function $h$ by
\begin{align}
   h := f_1^{-\frac{n-1}{2}}w .
\end{align}
Then, a direct computation shows that the function $h \bigl( r(p) \bigr)$ $( p \in N) $ satisfies the eigenvalue equation on $(N,g_N)$:
\begin{align*}
    - \Delta_{g_N} \bigl( h(r)  \bigr) 
    = -\left\{ \frac{\partial ^2}{\partial r^2}
    + (n-1) S(r) \frac{\partial }{\partial r} \right\} h(r)
    =  \left( \frac{(n-1)^2}{4} + 1 \right) h(r)
\end{align*}
and $h(r)\in L^2(N,dv_{g_N})$. 
Note that $dv_{g_N} = f_1^{\,n-1}(r)\,dr dv_{g_0}$, where $dv_{g_0}$ is the standard measure on the unit sphere $( S^{n-1}(1), g_0)$. 

Secondly, we shall construct a neighborhood of the origin of the desired manifold. 
Let $B_{{\bf R}^n}(0,r)$ be an open ball of radius $r$ and centered at the origin $0$ in the Euclidean space $({\bf R}^{n},g_{{\rm stand}})$ and denote by $\lambda _1 \bigl( B_{{\bf R}^n} (0,r) \bigr)$ the first Dirichlet eigenvalue of $B_{{\bf R}^n}(0,r)$. 
Since $\lim_{r\to +0} \lambda _1 \bigl( B_{{\bf R}^n} (0,r) \bigr) = \infty $ and $\lim_{r\to \infty } \lambda _1 \bigl( B_{{\bf R}^n} (0,r) \bigr) = 0$, there exists $ r_1 > 0 $ such that $ \lambda _1 \bigl( B_{{\bf R}^n} (0,r_1) \bigr) = \frac{(n-1)^2}{4} + 1$. 
Let $\widetilde{\varphi}_1$ be its associated positive-valued first eigenfunction. 
Since $\widetilde{\varphi}_1$ is a radial function, it can be written as $\widetilde{\varphi}_1 = H(r)$, where $r$ stands for the Euclidean distance to $0$.  
We note that $ H' < 0 $ on $(0,r_1]$. 

Thirdly, we shall connect two parts mentioned above; in view of $(54)$ and $(58)$, the function $h(t)$ also oscillates around $0$ and converges to $0$ as $t \to \infty$, and hence, there exist a constant $ r_2 > \max\{ r_1, 1\} $ such that $h(r_2) < 0 $ and $h'(r_2) < 0 $. 
Therefore, we can connect two functions, $H$ on $[0,r_1]$ and $h|_{[r_2,\infty)}$, by some function $\psi \in C^{\infty}[0,\infty)$ satisfying 
\begin{align}
  \psi (x)  = 
  \begin{cases}
    H(x)    & \qquad \mathrm{if}~~x \in [0 , r_1], \\
    h(x)    & \qquad \mathrm{if}~~x \in [r_2, \infty) ; 
  \end{cases}
\end{align}
and 
\begin{align}
    \psi ' (x) < 0 \qquad \mathrm{if}~~x \in [r_1 , r_2].
\end{align} 
Now, let us construct a function $f$ so that $\psi$ is an eigenfunction with eigenvalue $b_n := \frac{(n-1)^2}{4} + 1$ on $\bigl( {\bf R}^n, dr^2 + f^2(r) g_{S^{n-1}(1)} \bigr)$, that is,
\begin{align}
   \psi''(r) + (n-1) \frac{f'(r)}{f(r)}\psi'(r) = - b_n \psi(r).
\end{align}
By $(59)$ and $(60)$, we have $\psi ' < 0$ on $(0, r_2]$. 
Hence, we can solve the differential equation $(61)$ on the interval $[0,r_2]$ with the condition $f(r_1)=r_1$: 
\begin{align}
  f(r) = r_1 \exp \left\{ - \int_{r_1}^r 
      \frac{b_n\psi(s)+\psi''(s)}{(n-1)\psi'(s)} \,ds \right\}
     \qquad {\rm for}~~r \in [0,r_2].
\end{align}
Since $\psi |_{[0,r_1]}= H$ and $\widetilde{\varphi}_1 = H(r)$ on $B_{{\bf R}^n} (0,r)$, we see that $f(t) = t$ on $[0,r_1]$, and hence, $(B_{{\bf R}^n}(0,r_1), dr^2 + f^2(r)g_{S^{n-1}(1)})$ is a flat disk in $({\bf R}^{n},g_{{\rm stand}})$ with radius $r_1$; next, using this function $f$ on $[0,r_2]$, let us set 
\begin{align}
   f(r) = \frac{f(r_2) }{ f_1(r_2) } f_1 (r)
   \qquad {\rm for}~~r \in [r_2,\infty).
\end{align}
Then, this positive-valued function $f$ on $(0,\infty)$, defined by $(62)$ and $(63)$, satisfies the equation $(61)$, and hence, we see that $\psi$ is an eigenfunction with eigenvalue $b_n = \frac{(n-1)^2}{4} + 1$ on the manifold $(M,g):=\bigl( {\bf R}^n, dr^2 + f^2(r) g_{S^{n-1}(1)} \bigr)$; 
\begin{align}
    \frac{(n-1)^2}{4} + 1 \in \sigma_{{\rm p}}(-\Delta).
\end{align}
Moreover, from $(63)$, we have for $r \ge r_2$
\begin{align}
  & \nabla dr = S(r) \{ g - dr \otimes dr \} \\
  & K_{{\rm rad.}} = K(r) ,
\end{align}
and hence, Theorem $1.8$ $(1)$ follows from $(56)$ and $(65)$; Theorem $1.8$ $(3)$ follows from $(57)$ and $(66)$. 

In order to prove that $ - \Delta$ has no eigenvalue on the interval $\left( \frac{(n-1)^2}{4} ,\infty \right) $ except for the special number $ \frac{(n-1)^2}{4} + 1 $, we shall use the separation of variables: 
${\bf R}^n - \{ 0 \}$ is diffeomorphic to $(0,\infty) \times S^{n-1}(1)$ and 
we denote the eigenvalues of the Laplacian on the standard unit sphere $S^{n-1}(1)$ by
$$
   0=\lambda_0<\lambda_1\le \lambda_2\le \cdots
$$ 
with repetitions according to multiplicity; 
then, $-\Delta$ on $(M,g) = ({\bf R}^n,dr^2 + f^2(r) g_{S^{n-1}(1)})$ is unitarily equivalent to the infinite sum of the operators $-L_j$ on $L^2\bigl( (0,\infty),dx \bigr)$:
\begin{align*}
  &  - L_j = - \frac{d^2}{dx^2} + q_j \quad 
     \mathrm{on}~~L^2\bigl( (0,\infty),dx \bigr) \quad (j=0,1,2,\cdots);\\
  &  q_j(x) = \frac{(n-1)(n-3)}{4} \left( \frac{f'(x)}{f(x)} \right)^2 +
     \frac{n-1}{2}\frac{f''(x)}{f(x)} + \frac{\lambda_j}{f^2(x)}.
\end{align*}
Since 
\begin{align*}
   q_j(x) = \frac{(n-1)^2}{4} + \frac{k(n-1)\sqrt{(n-1)^2 + 4}}{2x}
            \sin (2x + c_n) + O(x^{-2}),
\end{align*}
Lemma $8.1$ implies that $-L_j$ has no eigenvalue on the interval $\left( \frac{(n-1)^2}{4} ,\infty \right) $ except for the special number $ \frac{(n-1)^2}{4} + 1 $. 
Thus we have proved Theorem $1.8$. 
\section{Remarks}
Tosio Kato \cite{Kato} proved that Schr\"odinger operator $ \Delta + V(x)$ on ${\bf R}^n$ has no eigenvalue $\lambda > K^2$ under the assumptions: $ V \in C^0({\bf R}^n)$, $\lim_{r \to \infty} V(x) = 0$ and $K=\lim_{r\to \infty} |rV(x)|$; 
Theorem $1.1$ seems to have a similar nature. 

In our theorems, we assume that there exists an open subset $U$ of $M$ with compact boundary $\partial U$ such that the outward pointing normal exponential map $\exp_{\partial U}:N^{+}(\partial U)\to M-\overline{U}$ induces a diffeomorphism. 
This condition is not essential. 
What matters is rather the existence of a function with special properties, such as $r$. 
The readers interested in this matter could pick up necessary conditions 
that should be satisfied by such a function from our proof above. 
We note that there are Donnelly's works (\cite{D4},\cite{D5}) from the viewpoint of an exhaustion function of $M$. 

In Section $8$, we have constructed a manifold with {\it one} end. 
However, we can also constructed a {\it two}-end-manifold satisfying the same properties; for that, it suffices to connect two copies of $(N,g_N)$ in Section $8$ by using a Riemannian product $[0,\frac{\pi}{\sqrt{b_n}}] \times S^{n-1}(\varepsilon)$ in a similar manner, where $S^{n-1}(\varepsilon) = \{ x \in {\bf R}^n \mid {\rm dist}(0,x) = \varepsilon \}$ and $\varepsilon > 0$ is a sufficiently small constant. 

\vspace{1mm}
\begin{flushleft}
Hironori Kumura\\ 
Department of Mathematics\\ 
Shizuoka University\\ 
Ohya, Shizuoka 422-8529\\ 
Japan\\
E-mail address: smhkumu@ipc.shizuoka.ac.jp
\end{flushleft}


\begin{thebibliography}{20}
\bibitem{A-K}K.~Akutagawa and H.~Kumura, 
{\it The uncertainty principle lemma under gravity}, 
arXiv:0812.4663, preprint.
\bibitem{A}F.~V.~Atkinson, 
{\it The asymptotic solution of second order differential equations}, 
Ann. Math. Pura. Appl., {\bf 37} (1954), 347--378.
\bibitem{A-U}M.~Arai and J.~Uchiyama, 
{\it On the von Neumann and Wigner potentials}, 
J. Differential Equations, {\bf 157} ((1999), 348--372.
\bibitem{D1}H.~Donnelly, 
{\it Eigenvalues embedded in the continuum for negatively curved manifolds}, 
Michigan Math. J., {\bf 28} (1981), 53--62.
\bibitem{D2}H.~Donnelly, 
{\it Negative curvature and embedded eigenvalues}, 
Math. Z., {\bf 203} (1990), 
301--308.
\bibitem{D3}H.~Donnelly, 
{\it Embedded eigenvalues for asymptotically flat surfaces}, 
Proceeding Symposia in Pure Mathematics, {\bf 54} (1993), Part 3, 169--177.
\bibitem{D4}H.~Donnelly, 
{\it Exhaustion functions and the spectrum of Riemannian manifolds}, 
Indiana Univ. Math. J. 46 (1997), 505--528.
\bibitem{D5}H.~Donnelly, 
{\it Spectrum of the Laplacian on asymptotically Euclidean spaces}, 
Michigan Math. J. 46 (1999), 101--111.
\bibitem{D-G}H.~Donnelly and N.~Garofalo, 
{\it Riemannian manifolds whose Laplacian have 
purely continuous spectrum}, Math. Ann., {\bf 293} (1992), 143--161.
\bibitem{E}D.~M.~Eidus, {\it The principle of limit amplitude}, 
Russian Math. Surveys, {\bf 24} (1969), no. 3, 97--167.
\bibitem{E}J.~Escobar, 
{\it On the spectrum of the Laplacian on complete Riemannian manifolds}, 
Comm. Partial Differential Equations, {\bf 11} (1986), 63--85.
\bibitem{E-F}J.~Escobar and A.~Freire, 
{\it The spectrum of the Laplacian of manifolds of positive curvature}, 
Duke Math. J., {\bf 65} (1992), 1--21.
\bibitem{G-W}R.~Green and H.~Wu, 
Function theory on manifolds which possess a pole, 
Lecture Note in Mathematics, vol.~699.
\bibitem{Karp}L.~Karp, 
{\it Noncompact manifolds with purely continuous spectrum}, 
Mich. Math. J., {\bf 31} (1984), 339--347.
\bibitem{Kato}T.~Kato, 
{\it Growth properties of solutions of the reduced wave equation with a variable coefficient}, 
Comm. Pure Appl. Math., {\bf 12} (1959), 403--426.
\bibitem{Kasue}A.~Kasue, 
Applications of Laplacian and Hessian comparison theorems, 
Geometry of geodesics and related topics (Tokyo, 1982), 333--386, Adv. Stud. Pure Math., 3.
\bibitem{K1}H.~Kumura, 
{\it On the essential spectrum of the Laplacian on complete manifolds}, 
J. Math. Soc. Japan, {\bf 49} (1997), 1--14.
\bibitem{K2}H.~Kumura, 
{\it A note on the absence of eigenvalues on negatively curved manifolds}, Kyushu J. Math., {\bf 56} (2002), 109--121.
\bibitem{M}K.~Mochizuki, 
{\it Growth properties of solutions of second order elliptic differential equations}, 
J. Math. Kyoto Univ., {\bf 16} (1976), 351--373.
\bibitem{P}M.~A.~Pinsky, 
{\it Spectrum of the Laplacian on a manifold of negative curvature II}, 
J. Differential Geometry, {\bf 14} (1979), 609--620.
\bibitem{R}S.~N.~Roze, 
{\it On the spectrum of an elliptic operator of second order}, 
Math. USSR. Sb., {\bf 9} (1969), 183--197.
\bibitem{T}T.~Tayoshi, 
{\it On the spectrum of the Laplace-Beltrami operator on a non-compact surface}, Proc. Japan. Acad., {\bf 47} (1971), 187--189.
\bibitem{W}S.~Wallach, 
{\it On the location of spectra of differential equations}, 
Amer. J. Math., {\bf 70} (1948), 833--841.
\end{thebibliography}
\end{document}